\documentclass[12pt]{amsart}
\usepackage{amssymb,times}

\makeatletter
\def\@strippedMR{}
\def\@scanforMR#1#2#3\endscan{
  \ifx#1M\ifx#2R\def\@strippedMR{#3}
  \else\def\@strippedMR{#1#2#3}
  \fi\fi}
\renewcommand\MR[1]{\relax\ifhmode\unskip\spacefactor3000 \space\fi
  \@scanforMR#1\endscan
  MR\MRhref{\@strippedMR}{\@strippedMR}}
\makeatother

\parskip=\medskipamount

\newtheorem*{Thm*}{Theorem}
\newtheorem{Thm}{Theorem}
\newtheorem{Cor}[Thm]{Corollary}
\newtheorem{Prop}[Thm]{Proposition}
\newtheorem{Lemma}[Thm]{Lemma}

\theoremstyle{definition}
\newtheorem{Defn}{Definition}

\newtheorem{Remark}{Remark}
\newtheorem{Example}{Example}

\newcommand{\mf}[1]{\mathbb{#1}}
\newcommand{\mc}[1]{\mathcal{#1}}
\newcommand{\mb}[1]{\mathbf{#1}}

\DeclareMathOperator{\Var}{\mathrm{Var}}
\DeclareMathOperator{\sgn}{\mathrm{sgn}}
\DeclareMathOperator{\supp}{\mathrm{supp}}
\DeclareMathOperator*{\esssup}{\mathrm{esssup}}

\newcommand{\norm}[1]{\left\Vert#1\right\Vert}
\newcommand{\abs}[1]{\left\vert#1\right\vert}
\newcommand{\chf}[1]{\mathbf{1}_{#1}}
\newcommand{\set}[1]{\left\{#1\right\}}
\newcommand{\ip}[2]{\left \langle #1, #2 \right \rangle}
\newcommand{\Exp}[1]{\mathbb{E} \left[#1\right]}
\newcommand{\Span}[1]{\mathrm{Span}\left\{#1\right\}}
\renewcommand{\phi}{\varphi}

\newcommand{\eps}{\varepsilon}
\newcommand{\Id}{\mathrm{Id}}
\newcommand{\subord}{\, \frame{$\ \vdash$} \,}

\allowdisplaybreaks[1]

\title{Generators of some non-commutative stochastic processes}
\author{Michael Anshelevich}
\thanks{This work was supported in part by NSF grants DMS-0900935 and DMS-1160849.}
\address{Department of Mathematics, Texas A\&M University, College Station, TX 77843-3368}
\email{manshel@math.tamu.edu}
\date{\today}

\begin{document}

\begin{abstract}
A fundamental result of Biane (1998) states that a process with freely independent increments has the Markov property, but that there are two kinds of free L{\'e}vy processes: the first kind has stationary increments, while the second kind has stationary transition operators. We show that a process of the first kind (with mean zero and finite variance) has the same transition operators as the free Brownian motion with appropriate initial conditions, while a process of the second kind has the same transition operators as a monotone L{\'e}vy process. We compute an explicit formula for the generators of these families of transition operators, in terms of singular integral operators, and prove that this formula holds on a fairly large domain. We also compute the generators for the $q$-Brownian motion, and for the two-state free Brownian motions.
\end{abstract}

\maketitle

\section{Introduction}

A L{\'e}vy process is a random process $\set{X(t): t \geq 0}$ whose increments $X(t) - X(s)$ are independent and stationary, in the sense that the distribution $\mu_{s,t}$ of $X(t) - X(s)$ depends only on $t-s$,
\[
\mu_{s,t} = \mu_{t-s}.
\]
A L{\'e}vy process is a Markov process, and its transition operators $\mc{K}_{s,t}$ defined via
\[
\Exp{f(X(t)) | s} = (\mc{K}_{s,t} f)(X(s))
\]
are also stationary, $\mc{K}_{s,t} = \mc{K}_{t-s}$; in fact
\[
\mc{K}_{s,t}(f)(x) = \int_{\mf{R}} f(x + y) \,d\mu_{t-s}(y).
\]
Here $\Exp{\cdot | s}$ is the conditional expectation onto time $s$. The maps $\set{\mc{K}_t : t \geq 0}$ form a semigroup, which has a generator $A$. In fact $A = \ell(- i \partial_x)$, where $\ell$ is the cumulant generating function of the process. It can also be expressed in terms of the L{\'e}vy measure of the process. See Section~\ref{Section:Classical} for more details.

In a groundbreaking paper \cite{BiaProcesses}, Biane showed that processes with freely independent increments, in the context of free probability \cite{VDN,Nica-Speicher-book}, are also Markov processes. He also noted that there are two distinct classes of such processes which can be called (additive) free L{\'e}vy process (Biane also investigated multiplicative processes, which we will not study here). Free L{\'e}vy processes of the first kind (FL1) have stationary increments, in the sense that each $X(t) - X(s)$ has distribution $\mu_{t-s}$. Then $\set{\mu_t: t \geq 0}$ form a semigroup with respect to free convolution $\boxplus$. These processes are Markov, but their transition operators typically are not stationary. Free L{\'e}vy processes of the second kind (FL2) have stationary transition operators $\mc{K}_t$, which form a semigroup, but their increments are typically not stationary: if $X(t) - X(s)$ has distribution $\mu_{s,t}$, then we only have the property $\mu_{s,t} \boxplus \mu_{t,u} = \mu_{s,u}$  (so that the measures form a free convolution hemigroup).

In this paper we compute the generators of free L{\'e}vy processes with finite variance. Since for FL1, the transition operators do not form a semigroup, they have a family of generators $\set{A_t: t \geq 0}$. In the case of FL2, there is a genuine generator $A$. If the distribution of the process has finite moments, using the free It{\^o} formula from \cite{AnsIto}, one can express the generators in terms of the $R$-transform, the free analog of $\ell$. However, it is unclear whether such an expression can be assigned a meaning in the absence of moments. But there is an alternative description. For a measure $\nu$, denote
\[
L_\nu(f)(x) = \int_{\mf{R}} \frac{f(x) - f(y)}{x - y} \,d\nu(y),
\]
a singular integral operator. A free convolution semigroup with finite variance is characterized by the free canonical pair $(\alpha, \rho)$, where $\alpha \in \mf{R}$, and (with appropriate normalization) $\rho$ is a probability measure. Further, denote by $\gamma_t$ the semicircular distribution at time $t$, so that $\rho \boxplus \gamma_t$ is the free analog of heat flow started at $\rho$. Then the generator of the corresponding free L{\'e}vy process of the first kind is
\begin{equation}
\label{Eq:Formula-Intro}
\alpha \partial_x + \partial_x L_{\rho \boxplus \gamma_t}.
\end{equation}
In fact, we show that for $\alpha = 0$, the full Markov structure of this process coincides with that of the free Brownian motion $\set{Y_t : t \geq 0}$ with $Y_0$ having distribution $\rho$. This statement clearly has no classical analogue.

In addition to free probability theory, there are only two other ``natural'' non-commutative probability theories \cite{Muraki-Natural-products}, the Boolean and the monotone. These theories do not, at least at this point, approach the wealth of structure of free or classical probability. However, one reason to study them is that they turn out to have unexpected connections to free probability. Indeed, we show that the generator of a free L{\'e}vy process of the second kind is $\alpha \partial_x + \partial_x L_{\rho}$, where $\rho$ now is the \emph{monotone} canonical measure. In fact, Biane in \cite{BiaProcesses} already noted that each FL2 is associated to a semigroup of analytic maps, and Franz in \cite{Franz-Monotone-Boolean} observed that exactly such semigroups are associated with monotone L{\'e}vy processes: the measures $\mu_{0,t}$ do not form a free semigroup, but they do form a semigroup under monotone convolution. In the monotone case itself there is no distinction between the L{\'e}vy processes of the first and second kind (so the free case is really special in this respect), and the generator of a monotone L{\'e}vy process is related to its monotone L{\'e}vy measure in the expected way \cite{Franz-Muraki-Markov-monotone}.

We also compute the generators of the $q$-Brownian motion. This non-commutative process was constructed in \cite{BozSpeBM1}, and investigated in detail in \cite{BKSQGauss}. Building on the work of \cite{Donati-Martin}, we prove a functional It{\^o} formula for it (for polynomial functions), from which the formulas for generators easily follow.

We note that the study of ``time-dependent generators'', or more usually the inverse problem---how to reconstruct $\set{\mc{K}_{s,t}}$ from $\set{A_t}$---goes back to \cite{Kato-Non-autonomous}. This is typically formulated at the linear non-autonomous Cauchy problem, and a significant amount of general results on its solution is known, see for example Section VI.9 of \cite{Engel-Nagel-book}, \cite{Neidhardt-Zagrebnov-non-autonomous}, and their extensive references. We do not use these general results in the paper, but this may be a matter for further study.

The paper is organized as follows. After the introduction and some general results in Section~\ref{Section:Preliminaries}, in Section~\ref{Section:Classical} we give a short overview of the generators for classical processes. The next section, covering free L{\'e}vy processes, is the main part of the paper. We show that transition operators for such a process form a strongly continuous contractive family on $L^p(\mf{R}, \,dx)$, and that their generators are given by formula~\eqref{Eq:Formula-Intro} on large domains in $L^p(\mf{R}, \,dx)$ and $C_0(\mf{R})$. In Section~\ref{Subsec:Semigroups} we find the closures of these generators. We also compute the generators of FL2 processes. In section~\ref{Subsec:Further}, we show that the operator $L_\nu$ itself is an isometry between certain $L^2$ spaces, and compute the ``carr{\'e} du champ'' operator corresponding to $\partial_x L_\nu$. In Section~\ref{Section:q}, we compute the It{\^o} formula and generators for the $q$-Brownian motion, and in a short final section we apply similar analysis to the two-state free Brownian motions from \cite{Ans-Two-Brownian}.

\textbf{Acknowledgements.} This work was prompted by a discussion with W{\l}odek Bryc about his paper \cite{Bryc-Markov-Meixner}; I am grateful to W{\l}odek for showing me an early version of that paper. I have discussed various aspects of this article with a large number of people. Thanks to Dominique Bakry, Todd Kemp, Michel Ledoux, Conni Liaw, Alex Poltoratski, and Sergei Treil for helpful comments. Thanks especially to J.C. Wang and the referee for a correction in Remark~\ref{Remark:Definition}. Finally, I am grateful to the Erwin Schr{\"o}dinger Institute, and to the Universit{\'e} Paul Sabatier, where part of this work was completed.

\section{Preliminaries}
\label{Section:Preliminaries}

\subsection{Generalities and definitions}
\label{Subsec:Generalities}

Let $(\mc{M}, \mf{E})$ be a tracial non-commutative probability space, where $\mc{M}$ is a von Neumann algebra and $\mf{E}$ is a tracial normal state on it. Possibly unbounded random variables are self-adjoint elements of the algebra $\widetilde{\mc{M}}$ of operators affiliated to $\mc{M}$.

A \emph{process} is a family of (possibly non-commutative) random variables $\set{X(t): t \geq 0}$ in a (possibly non-commutative) probability space $(\mc{M}, \mf{E})$.

We will assume that $X(0) = 0$, and will denote by $\mu_{s,t}$ the distribution of $X(t) - X(s)$ with respect to $\mf{E}$ (for $s \leq t$), $\mu_t = \mu_{0,t}$ the distribution $X(t)$, and $\mu = \mu_1$. If $\star$ is a convolution operation corresponding to some non-commutative independence, and the increments of $\set{X(t)}$ are independent in that sense, then
\[
\mu_{s,t} \star \mu_{t,u} = \mu_{s,u}.
\]

For an unbounded operator $X$, we will denote by $\mc{D}(X)$ its domain, and by $(X, \mc{D})$ its restriction to a smaller domain $\mc{D}$.

\begin{Defn}
\label{Defn:Generator}
For a family of distributions $\set{\mu_t}$, we say that the functional $L_t$ is \emph{its generator at time $t$ with domain $\mc{D}(L_t)$} if
\[
\partial_t \int_{\mf{R}} f(x) \,d\mu_t(x) = L_t[f]
\]
for $f \in \mc{D}(L_t)$. Frequently,
\[
L_t[f] = \int_{\mf{R}} (A_t f)(x) \,d\mu_t(x)
\]
for an operator $A_t$. If $\set{X(t)}$ is a process with distributions $\set{\mu_t}$, this is equivalent to
\begin{equation}
\label{Eq:Generator-expectations}
\partial_t \Exp{f(X(t))} = \Exp{(A_t f)(X(t))}.
\end{equation}
Note however that this property does not determine $A_t$, even on $\mc{D}(L_t)$.

For operators $\set{\mc{K}_{s,t}}$ on a Banach space $\mc{A}$, we write
\[
\left. \frac{\partial}{\partial t} \right|_{t=s} \mc{K}_{s,t} = A_s
\]
if
\begin{equation}
\label{Eq:Generator}
\lim_{h \rightarrow 0^+} \norm{\frac{1}{h} \left( \mc{K}_{s, s+h} f - \mc{K}_{s,s} f \right) - A_s f} = 0.
\end{equation}
In this case we say that $A_s$ is the generator of the family $\set{\mc{K}_{s,t}}$  at time $s$. Its domain $\mc{D}(A_s) \subset \mc{A}$ consists of all $f \in \mc{A}$ for which the limit~\eqref{Eq:Generator} holds.

Now suppose that the process $\set{X(t)}$ is a Markov process. That is, denoting $\Exp{\cdot | \leq s}$ the $\mf{E}$-preserving conditional expectation onto the von Neumann algebra generated by $\set{X_u : u \leq s}$, for any $f \in L^\infty(\mf{R}, \,dx)$, $\Exp{f(X(t)) | \leq s}$ is in the von Neumann algebra generated by $X(s)$. (See the Introduction and Section~4 of \cite{BiaProcesses} for more details, and also for a weaker requirement, sufficient for our purposes, that the classical version of $\set{X(t)}$ is a Markov process.) In this case, the corresponding transition operators are determined by
\[
\Exp{f(X(t)) | \leq s} = (\mc{K}_{s,t} f)(X(s))
\]
We say that the operator $A_s$ is the generator of the process at time $s$ if it is the generator of its family of transition operators. Note that if $A_t$ exists, it has the property in equation~\eqref{Eq:Generator-expectations}.
\end{Defn}

\begin{Prop}
\label{Prop:Basic}
Let $(\mc{A}, \norm{\cdot})$ be a Banach space, $\set{\mc{K}_{s,t}}$ a family of contractions on $\mc{A}$ such that
\[
\mc{K}_{s,t} \circ \mc{K}_{t,v} = \mc{K}_{s,v}, \qquad \mc{K}_{s,s} = I,
\]
and $\mc{K}_{s,t}$ is strongly continuous in $t$. Let $\set{A_t}$ be the generators of $\set{\mc{K}_{s,t}}$ in the sense of equation~\eqref{Eq:Generator}, and consider a subspace $\mc{D} \subset \bigcap_{t} \mc{D}(A_t)$ such that for any $f \in \mc{D}$, $A_t f$ is a continuous function of $t$.
\begin{enumerate}
\item
Each $A_t$ is dissipative and closable.
\item
Let $\mc{B} \subset \mc{A}$ be another subspace such that $\mc{D} \subset \mc{B}$, and $\norm{\cdot}_{\mc{B}}$ be another norm on $\mc{B}$ such that $\mc{D}$ is $\norm{\cdot}_{\mc{B}}$-dense in $\mc{B}$, $\norm{f} \leq \norm{f}_{\mc{B}}$, and for $f \in \mc{D}$,
\begin{equation}
\label{Eq:B-bounded}
\norm{A_t f} \leq \norm{f}_{\mc{B}}.
\end{equation}
Then equation~\eqref{Eq:Generator} holds for $f \in \mc{B}$, so that $\mc{B} \subset \mc{D}(A_t)$ for all $t$.
\item
The closure $\overline{\mc{D}}^{\norm{\cdot}_A}$ of $\mc{D}$ in the sup-graph norm
\[
\norm{f}_A = \norm{f} + \sup_t \norm{A_t f}
\]
is in $\mc{D}(A_t)$ for all $t$.
\end{enumerate}
\end{Prop}

\begin{Remark}
Note that strong continuity of $\mc{K}_{s,t}$ does not imply continuity of $\set{A_t}$. Indeed, already in one dimension, if $\mc{K}_{s,t} = e^{f(t) - f(s)}$, then $A_t = f'(t)$.
\end{Remark}

\begin{proof}
For part (a), recall from Section X.8 of \cite{ReeSim2} that for $f \in \mc{A}$, a normalized tangent functional $\phi_f$ is an element of $\mc{A}^\ast$ such that $\norm{\phi_f} = \norm{f}$ and $\phi_f[f] = \norm{f}^2$. For any such functional,
\[
\begin{split}
\Re \phi_f[A_s f]
& = \lim_{h \rightarrow 0^+} \frac{1}{h} \Re \phi_f[\mc{K}_{s, s+h} f - f]
\leq \lim_{h \rightarrow 0^+} \frac{1}{h} \left( \abs{\phi_f[\mc{K}_{s, s+h} f]} - \norm{f}^2 \right) \\
& \leq \lim_{h \rightarrow 0^+} \frac{1}{h} \left( \norm{f} \cdot \norm{\mc{K}_{s, s+h} f} - \norm{f}^2 \right)
\leq 0
\end{split}
\]
since $\mc{K}_{s, s+h}$ is a contraction, so $A_s$ is dissipative. Combining this with Theorem~II.3.23 and Proposition~II.3.14(iv) of \cite{Engel-Nagel-book}, it follows that $A_s$ is closable.

For part (b), we first note that for $s \leq t$, since $\mc{K}_{s,t}$ is a contraction, for $f \in \mc{D}$,
\begin{equation*}
\begin{split}
\lim_{h \rightarrow 0^+} \norm{\frac{1}{h} \left( \mc{K}_{s, t+h} f - \mc{K}_{s,t} f \right) - \mc{K}_{s,t}(A_t f)}
& = \lim_{h \rightarrow 0^+} \norm{\mc{K}_{s,t} \left( \frac{1}{h} \left( \mc{K}_{t, t+h} f - f \right) - A_t f\right)} \\
& \leq \lim_{h \rightarrow 0^+} \norm{\frac{1}{h} \left( \mc{K}_{t, t+h} f - f \right) - A_t f}
= 0,
\end{split}
\end{equation*}
so
\[
\partial_t \mc{K}_{s,t}(f) = \mc{K}_{s,t}(A_t f).
\]
Also, since $\mc{K}_{s,v}$ is a contraction, $\mc{K}_{s,v} A_v f$ is continuous in $v$. Therefore we have the Riemann integral identity
\[
\mc{K}_{s,t}(f) = f + \int_s^t \mc{K}_{s,v}(A_v f) \,dv.
\]

Since for $f \in \mc{D}$, \eqref{Eq:B-bounded} holds, $A_t$ has a continuous extension $(\mc{B}, \norm{\cdot}_{\mc{B}}) \rightarrow \mc{A}$ satisfying the same property. We will show that this continuous extension $\tilde{A}_t$ coincides with $A_t$.

Fix $g \in \mc{B}$ and a time $t$. For each $\eps > 0$, we can find a $f \in \mc{D}$ such that $\norm{f - g}_{\mc{B}} < \eps$, so that $\norm{f - g} < \eps$,
\[
\norm{\mc{K}_{s,t} f - \mc{K}_{s,t} g} < \eps,
\]
and
\[
\norm{\mc{K}_{s,v} (\tilde{A}_v  f) - \mc{K}_{s,v} (\tilde{A}_v g)} \leq \norm{\tilde{A}_v f - \tilde{A}_v g} < \eps
\]
for all $s \leq v \leq t$. Then
\[
\norm{\mc{K}_{s,t} g - g - \int_s^t \mc{K}_{s,v}(\tilde{A}_v g) \,dv} < 2 \eps + (t-s) \eps.
\]
So
\begin{equation}
\label{Eq:Riemann-integral}
\mc{K}_{s,t}(g) = g + \int_s^t \mc{K}_{s,v}(\tilde{A}_v g) \,dv
\end{equation}
(in particular, the integral is well defined), and
\[
A_t g = \left. \frac{\partial}{\partial t} \right|_{t=s} \mc{K}_{s,t} g = \tilde{A}_t g.
\]

Finally, for part (c) we take $\mc{B} = \overline{\mc{D}}^{\norm{\cdot}_A}$ and $\norm{\cdot}_{\mc{B}} = \norm{\cdot}_A$, and apply part (b).
\end{proof}

\begin{Remark}
Under the assumptions of the preceding proposition, from equation~\eqref{Eq:Riemann-integral} we also get
\[
\left. \frac{\partial}{\partial s} \right|_{s=t} \mc{K}_{s,t} g = - A_s g,
\]
\[
A_s \mc{K}_{s,t}(g) = - \partial_s \mc{K}_{s,t}(g),
\]
which in turn implies
\[
\mc{K}_{s,t}(g) = g + \int_s^t A_v \mc{K}_{v,t}(g) \,dv.
\]
\end{Remark}

\begin{Lemma}
\label{Lemma:Generator-martingale}
Let $\set{X(t)}$ be a non-commutative Markov process, with transition operators $\set{\mc{K}_{s,t}}$ and distributions $\set{\mu_t}$. Let $(\mc{A}, \norm{\cdot})$ be either $(C_0(\mf{R}), \norm{\cdot}_\infty)$ or $L^p(\mf{R}, \,dx)$, and suppose that $\set{\mc{K}_{s,t}}$, their generators $\set{A_s}$, and $\mc{D} \subset \mc{A}$ satisfy the hypotheses of Proposition~\ref{Prop:Basic}. Then for $f \in \mc{D}$,
\begin{equation}
\label{Eq:Martingale}
f(X(t)) - \int_0^t (A_v f)(X(v)) \,dv
\end{equation}
is a martingale. Conversely, suppose $\set{B_s}$ is another family of operators strongly continuous on $\mc{D}$ such that \eqref{Eq:Martingale} is a martingale. Then for $f \in \mc{D}$, $A_s f = B_s f$ in the restriction of $\mc{A}$ to $\supp(\mu_s)$.
\end{Lemma}

\begin{proof}
As shown in Proposition~\ref{Prop:Basic}, under these assumptions, for $f \in \mc{D}$
\[
\mc{K}_{s,t}(f) = f + \int_s^t \mc{K}_{s,v}(A_v f) \,dv.
\]
It follows that
\[
\begin{split}
& \Exp{ \left. f(X(t)) - \int_0^t (A_v f)(X(v)) \,dv \right| \leq s} \\
&\qquad = (\mc{K}_{s,t}(f))(X(s)) - \int_0^s (A_v f)(X(v)) \,dv - \int_s^t \mc{K}_{s,v}(A_v f)(X(s)) \,dv \\
&\qquad = f(X(s)) - \int_0^s (A_v f)(X(v)) \,dv
\end{split}
\]
and the process is a martingale.

Conversely, suppose that $f(X(t)) - \int_0^t (B_v f)(X(s)) \,dv$ is a martingale. The last equality then implies that
\[
\mc{K}_{s,t}(f)(X(s)) = f(X(s)) + \int_s^t \mc{K}_{s,v}(B_v f)(X(s)) \,dv.
\]
It follows that in $C_0(\supp(\mu_s))$ or $L^p(\supp(\mu_s), \,dx)$,
\[
\mc{K}_{s,t}(f) = f + \int_s^t \mc{K}_{s,v}(B_v f) \,dv,
\]
and therefore in this space
\[
B_v f = \left. \frac{\partial}{\partial t} \right|_{t=v} \mc{K}_{v,t}(f) = A_v f.
\]
\end{proof}

\subsection{Cumulants}
\label{Subsec:Cumulants}

Let $\set{\mu_t}$ be a convolution semigroup with respect to some convolution operation $\star$. In all cases we will consider, $\mu_0 = \delta_0$, $\mu_t[x] = t \cdot \mu[x]$, and $\set{\mu_t}$ is weakly continuous. Almost by definition (see Property (K1') in Section~3 of \cite{Hasebe-Saigo-Monotone-cumulants}), the cumulant functional of $\mu$ corresponding to the convolution operation $\star$ is
\begin{equation}
\label{Eq:Cumulant-definition}
C_\mu[f] = \left. \frac{\partial}{\partial t} \right|_{t=0} \mu_t[f].
\end{equation}
This approach works for all the convolutions associated to natural types of independence (tensor, free, Boolean, monotone), but also for other operations such as the $q$-convolution from \cite{AnsQCum}.

\begin{Prop}
Assume that $\mu$ has finite moments of all orders. Then, at least on the space $\mc{P}$ of polynomials,
\begin{equation}
\label{Eq:Cumulant}
C_\mu[f] = \alpha f'(0) + \int_{\mf{R}} \frac{f(y) - f(0) - y f'(0)}{y^2} \,d\rho(y)
\end{equation}
for a finite measure $\rho$.
\end{Prop}

\begin{proof}
First note that $C_\mu[1] = 0$ and $C_\mu[x] = \mu[x] = \alpha$ for some $\alpha$. Since each $\mu_t$ is positive, it follows that $C_\mu$ is a conditionally positive function on polynomials, so it has the canonical representation $C_\mu[x^n] = \rho[x^{n-2}]$ for $n \geq 2$, where $\rho$ is a finite measure, the canonical measure of the semigroup. We compute, for $f(x) = x^n$
\[
C_\mu[x^n] = \rho[x^{n-2}] = \alpha f'(0) + \rho \left[ \frac{f(x) - f(0) - x f'(0)}{x^2} \right] = \alpha f'(0) +  \rho \left[ \left. \frac{\partial}{\partial x} \right|_{x=0} \frac{f(x) - f(y)}{x - y} \right]
\]
and this formula is also valid for $f(x) = 1$ and $x$.
\end{proof}

\section{Classical L{\'e}vy processes.}
\label{Section:Classical}

See Section~1.2 of \cite{Bertoin-book} (except for a small misprint) for the following results.

\begin{Thm*}
Let $\set{X(t)}$ be a L{\'e}vy process corresponding to the convolution semigroup $\set{\mu_t}$. Denote
\[
\ell(\theta) = \log \Exp{e^{i \theta X}} = \log \int_{\mf{R}} e^{i \theta x} \,d\mu(x)
\]
the cumulant generating function of the process. Then the generator of the process is the pseudo-differential operator $\ell(- i \partial_x)$ with dense domain
\[
\set{f \in L^2(\mf{R}, dx): \int_{\mf{R}} \abs{\ell(\theta)}^2 \abs{\hat{f}(\theta)}^2 \,d\theta < \infty}.
\]
In other words, if the process has the L{\'e}vy-Khintchine representation
\[
\ell(\theta) = i \alpha \theta - \frac{1}{2} V \theta^2 + \int_{\mf{R} \backslash \set{0}} (e^{i y \theta} - 1 - i y \theta \chf{\abs{y} < 1}) \Pi(dy),
\]
then
\[
A f(x) = \alpha f'(x) + \frac{1}{2} V f''(x) + \int_{\mf{R}} \Bigl(f(x+y) - f(x) - \chf{\abs{y} < 1} y f'(x) \Bigr) \Pi(dy).
\]
\end{Thm*}

If $\mu$ has mean $\alpha$ and finite variance, we also have the Kolmogorov representation,
\[
\ell(\theta) = i \alpha \theta + \int_{\mf{R}} (e^{i y \theta} - 1 - i y \theta ) y^{-2} \,d\rho(y),
\]
where $\rho$ is the canonical measure. In this case the generator is
\begin{equation}
\label{Eq:Generator-classical}
A f (x) = \alpha f'(x) + \int_{\mf{R}} \frac{f(x+y) - f(x) - y f'(x)}{y^2} \,d\rho(y).
\end{equation}
If the process has (say) finite exponential moments, we have moveover
\[
\ell(\theta) = \sum_{n=1}^\infty \frac{1}{n!} c_n (i \theta)^n,
\]
where $\set{c_n}$ are the cumulants \cite{Shi} of (the distribution of) the process. So the generator of the process is
\[
\sum_{n=1}^\infty \frac{1}{n!} c_n \partial_x^n.
\]
Note also that $c_1 = \alpha$ and
\[
c_n = \int_{\mf{R}} x^{n-2} \,d\rho(x)
\]
for $n \geq 2$.

\section{Free L{\'e}vy processes}
\label{Section:Free}

\subsection{Background}
\label{Subsec:Background}

Let $\mu$ be a probability measure on $\mf{R}$. Its Cauchy transform is the analytic function $G_\mu:\mf{C}^+ \rightarrow \mf{C}^-$ defined by
\[
G_\mu(z) = \int_{\mf{R}} \frac{1}{z - x} \,d\mu(x).
\]

We will also denote $F_\mu(z) = \frac{1}{G_\mu(z)}$, so that $F_\mu : \mf{C}^+ \rightarrow \mf{C}^+$.

$G_\mu$ is invertible in a Stolz angle near infinity, and Voiculescu's $R$-transform is defined by
\[
R_\mu(z) = G_\mu^{-1}(z) - \frac{1}{z}.
\]
A free convolution $\mu_1 \boxplus \mu_2$ of two probability measures $\mu_1, \mu_2$ is characterized by the property that
\[
G_{\mu_1 \boxplus \mu_2}^{-1}(z) = G_{\mu_1}^{-1}(z) + R_{\mu_2}(z).
\]
$\mu$ is $\boxplus$-infinitely divisible if and only it can be included as $\mu = \mu_1$ in a free convolution semigroup $\set{\mu_t : t \geq 0}$, $\mu_s \boxplus \mu_t = \mu_{s+t}$. This is the case if and only if $R_\mu$ extends to an analytic function $R_\mu: \mf{C}^+ \rightarrow \mf{C}^+ \cup \mf{R}$. In this case, we have the free L{\'e}vy-Khintchine representation (Theorem~5.10 of \cite{BV93})
\[
R_\mu(z) = \alpha + \int_{\mf{R}} \frac{z + x}{1 - x z} \,d\nu(x).
\]
Moreover, if $\mu$ has finite variance, we also have the free Kolmogorov representation,
\[
R_\mu(z) = \alpha + \int_{\mf{R}} \frac{z}{1 - x z} \,d\rho(x).
\]
Here $\alpha$ is the mean of $\mu$, $\rho$ is a finite measure, the free canonical measure for the semigroup $\set{\mu_t}$, and $(\alpha, \rho)$ is the free canonical pair. For convenience, throughout most of the paper we will rescale time so that the variance
\begin{equation}
\label{Eq:Normalization}
\Var[\mu = \mu_1] = 1
\end{equation}
in which case $\rho$ is a probability measure.

We will also encounter two other convolution operations, the monotone convolution $\triangleright$ and the Boolean convolution $\uplus$, determined by
\[
F_{\mu_1 \triangleright \mu_2}(z) = F_{\mu_1}(F_{\mu_2}(z))
\]
and
\[
F_{\mu_1 \uplus \mu_2}(z) = F_{\mu_1}(z) + F_{\mu_2}(z) - z.
\]

We will denote
\[
d\gamma_t(x) = \frac{1}{2 \pi t} \sqrt{4 t - x^2} \chf{[-2 \sqrt{t}, 2 \sqrt{t}]}(x) \,dx
\]
the semicircular distributions, the analogs of the normal distributions in free probability. They form a free convolution semigroup with the free canonical pair $(0, \delta_0)$.

\begin{Remark}
\label{Remark:Stripping}
If $\nu$ is a probability measure, there exists a probability measure $\mu = \Phi_t[\nu]$ with mean zero and variance $t$ such that
\[
F_{\Phi_t[\nu]}(z) = z - t G_\nu(z).
\]
In particular, for the map $\Phi = \Phi_1$, see \cite{Belinschi-Nica-Free-BM} and Proposition~2.2 of \cite{Maa92}, and
\[
\Phi_t[\nu] = \Phi[\nu]^{\uplus t}.
\]
Conversely, if $\mu$ is a probability measure with mean $\alpha$ and variance $\beta > 0$, there exists a probability measure $\nu = \mc{J}[\mu]$ such that
\begin{equation}
\label{Eq:J-definition}
F_\mu(z) = z - \alpha - \beta G_{\mc{J}[\mu]}(z).
\end{equation}
Note that $\mc{J} \circ \Phi_t = \Id$, while $\Phi_t[\mc{J}[\mu]] = \mu$ if $\mu$ has mean zero and variance $t$. If $\nu$ has finite moments of all orders, its Cauchy transform has a continued fraction expansion
\[
G_\mu(z) =
\cfrac{1}{
z - \alpha_0 - \cfrac{\beta_1}{
z - \alpha_1 - \cfrac{\beta_2}{
z - \alpha_2 - \cfrac{\beta_3}{
z - \ldots}}}}.
\]
Here $\beta_0 = 1$, $\alpha_0$ is the mean of $\mu$, $\beta_1$ is the variance of $\mu$, and in general $\set{\alpha_n, \beta_n}$ are its Jacobi parameters. Then for $\nu = \Phi_t[\mu]$,
\[
G_\nu(z) =
\cfrac{1}{
z - \cfrac{t}{
z - \alpha_0 - \cfrac{\beta_1}{
z - \alpha_1 - \cfrac{\beta_2}{
z - \alpha_2 - \cfrac{\beta_3}{
z - \ldots}}}}},
\]
while $\mc{J}$ is the inverse map, namely coefficient stripping \cite{Damanik-Simon-periodic}. Note that there are also related maps which involve finite rather than only probability measures, but because of the normalization~\eqref{Eq:Normalization}, we will not need to consider them.
\end{Remark}

\subsection{Transition operators}
\label{Subsec:Transition}

The following is a reformulation of Theorem~3.1 of \cite{BiaProcesses}.

\begin{Thm*}
Let $X$ and $Y$ be freely independent. Then the transition operator $\mc{K}$ defined via
\[
\Exp{f(X+Y) | X} = (\mc{K} f)(X)
\]
is a map on $C_0(\mf{R})$ (which extends to a map on $L^\infty(\mf{R}, \,dx)$) such that for any $z \in \mf{C} \setminus \mf{R}$
\[
\mc{K} \left[ \frac{1}{z - x} \right]
= \frac{1}{F(z) - x}
= \frac{1}{F_\nu(z) - x}.
\]
Here $F(z) = F_\nu(z)$ for some probability measure $\nu$, and $F$ is uniquely determined by
\[
G_{X+Y}(z) = G_X(F(z)).
\]
\end{Thm*}

\begin{Prop}
For processes with freely independent increments, the transition operators $\mc{K}$ have the form
\[
\mc{K}[f](x) = \int_{\mf{R}} f(y) \,d(\delta_x \triangleright \nu)(y).
\]
For the free L{\'e}vy processes of the second kind, for $\mc{K}_t$ we have $\nu_t = \mu_t$. For the free L{\'e}vy processes of the first kind, for $\mc{K}_{s,t}$, $\nu_{s,t}$ is determined by
\[
G_{\mu_t}(z) = G_{\mu_s}(F_{\nu_{s,t}}(z)).
\]
In other words, $\nu = \nu_{s,t} = \mu_{t-s} \subord \mu_s$, the subordination distribution, see \cite{Lenczewski-Decompositions-convolution,Nica-Subordination}.
\end{Prop}

\begin{proof}
According to Theorem~3.1 of \cite{BiaProcesses},
\[
\mc{K} \left[ \frac{1}{z - x} \right]
= \frac{1}{F(z) - x}
= \frac{1}{F_\nu(z) - x}
= \frac{1}{F_{\delta_x \triangleright \nu}(z)}
= G_{\delta_x \triangleright \nu}(z)
= \int_{\mf{R}} \frac{1}{z - y} \,d(\delta_x \triangleright \nu)(y),
\]
and, still according to \cite{BiaProcesses}, this property entirely determines $\mc{K}$. For FL2, $F_{s,t} = F_{t-s}$ and $F_{\mu_t} = F_{\mu_s} \circ F_{t-s}$, so for $s = 0$,
\[
F_{\mu_t} = F_{\delta_0} \circ F_t = F_t = F_{\nu_t}.
\]
For FL1, $\nu_{s,t} = \mu_{t-s} \subord \mu_s$ by definition.
\end{proof}

\begin{Remark}
Note that
\[
\mc{K}(x, dy) = (\delta_x \triangleright \nu)(y) = (\delta_x \uplus \nu)(y) = (\nu \uplus \delta_x)(y).
\]
In fact, measures $\delta_x \triangleright \nu$ are well-known in classical spectral theory. Indeed, if $X$ is an operator with cyclic vector $\xi$ and corresponding distribution $\nu$, then $\delta_x \triangleright \nu$ is the distribution with respect to $\xi$ of the rank-one perturbation $X - x \ip{\cdot}{\xi} \xi$. Finally, note that for the classical processes and convolution, we can also write
\[
\mc{K}[f](x)= \int_{\mf{R}} f(y) \,d\mu(y - x) = \int_{\mf{R}} f(y) \,d (\delta_x \ast \mu)(y)
\]
\end{Remark}

\begin{Prop}
\label{Prop:Contraction}
$\mc{K}$ is a contraction on each $L^p(\mf{R}, dx)$ for $1 \leq p \leq \infty$.
\end{Prop}

\begin{proof}
For $f \in L^\infty(\mf{R}, dx)$,
\[
\norm{\mc{K}[f]}_\infty
= \esssup_{x \in \mf{R}} \int_{\mf{R}} \abs{f(y)} \,d(\delta_x \triangleright \nu)(y)
\leq \esssup_{x \in \mf{R}} \norm{f}_\infty = \norm{f}_\infty,
\]
so $\mc{K}$ is a contraction on $L^\infty(\mf{R}, dx)$. On the other hand, Alexandrov's averaging theorem (Theorem~11.8 from \cite{Simon-Trace-ideals}) states that for $f \in L^1(\mf{R}, \,dx)$
\[
\int_{\mf{R}} \mc{K}[f](x) \,dx = \int_{\mf{R}} \left( \int_{\mf{R}} f(y) \,d(\delta_x \triangleright \nu)(y) \right) \,dx = \int_{\mf{R}} f(x) \,dx,
\]
so that $\norm{\mc{K} f}_1 = \norm{f}_1$ for $f \geq 0$, and $\mc{K}$ is a contraction on $L^1(\mf{R}, \,dx)$. For $1 < p < \infty$, the result now follows by Riesz-Thorin interpolation, see Section~IX.4 from \cite{ReeSim2}.
\end{proof}

\begin{Remark}
\label{Remark:Definition}
Unless stated otherwise, we will work in the Banach space $C_0(\mf{R})$ of continuous functions converging to zero at infinity, with the maximum norm. Also, again until stated otherwise, we will denote
\begin{equation}
\label{Eq:Domain-rational-functions}
\mc{D} = \Span{ \frac{1}{z_1 - x} - \frac{1}{z_2 - x} : z_1, z_2 \in \mf{C} \setminus \mf{R}}.
\end{equation}
We will also denote by $C_c(\mf{R})$ the compactly supported continuous functions, and by $\mc{P}$ the polynomials.

For $1 \leq p \leq \infty$, denote
\[
\norm{f}_{k, p} = \sum_{i=0}^k \norm{f^{(i)}}_p \sim \norm{f^{(k)}}_p + \norm{f}_p
\]
the Sobolev norm. For $p < \infty$, denote by $W^{k, p}(\mf{R})$ the corresponding Sobolev space, while for $p = \infty$ we will denote by $W^{k, \infty}$ the corresponding subspace of $C_0(\mf{R})$. Note that we will identify the Lipschitz norm
\[
\sup_{x \neq y} \abs{\frac{f(x) - f(y)}{x - y}}
\]
with  $\norm{f'}_\infty$, since a Lipschitz function is differentiable almost everywhere.

Finally, we abbreviate
\[
W^\infty = \set{f \in C_0(\mf{R}) | f'' \in C_0(\mf{R})}
\]
with norm $\norm{f}_\infty + \norm{f''}_\infty$, and
\[
W^p = W^\infty \cap W^{1,p}
\]
with norm
\begin{equation}
\label{Eq:Norm-four-terms}
\norm{f}' = \norm{f}_\infty + \norm{f''}_\infty + \norm{f}_p + \norm{f'}_p.
\end{equation}
\end{Remark}

The following argument is reminiscent of Section~2.3 from \cite{Liaw-Treil-Rank-one}.

\begin{Lemma}
\label{Lemma:W^p-density}
$\mc{D}$ is dense in $W^p$, $1 \leq p \leq \infty$, $L^p(\mf{R}, \,dx)$, $1 \leq p < \infty$, and $C_0(\mf{R})$, with their respective norms.
\end{Lemma}

\begin{proof}
We will prove that $\mc{D}$ is dense in $W^p$; the other arguments are similar, and more standard. Note that $\mc{D} \subset W^p$.

\textbf{Step 1.}
By an elementary ``cut-off plus smoothing'' argument,
\begin{equation}
\label{Eq:Compactly-supported}
C_c(\mf{R}) \cap W^p = C_c(\mf{R}) \cap W^\infty
\end{equation}
is dense in $W^p$ with respect to the norm~\eqref{Eq:Norm-four-terms}. So it suffices to check that every element of this space can be approximated by elements of $\mc{D}$.

\textbf{Step 2.}
For
\[
P_\eps(x) = \frac{1}{\pi} \frac{\eps}{x^2 + \eps^2}
\]
the Poisson kernel and $f$ in the set~\eqref{Eq:Compactly-supported}, we know that
\begin{equation}
\label{Eq:Poisson-derivative}
(P_\eps \ast f)'' = P_\eps \ast f'' = P_\eps'' \ast f.
\end{equation}
Moreover, as $\eps \rightarrow 0^+$, $(P_\eps \ast f) \rightarrow f$ and $(P_\eps \ast f)'' \rightarrow f''$ uniformly, and so (since the support of $f$ is compact) also in norm~\eqref{Eq:Norm-four-terms}.

\textbf{Step 3.}
For a fixed $\eps$, if
\[
\sum \frac{1}{\pi} \frac{\eps}{(x - a_i) + \eps^2} f(a_i) \Delta_i
\]
is a Riemann sum for $(P_\eps \ast f)$, then by \eqref{Eq:Poisson-derivative}
\[
\sum \frac{1}{\pi} \left(\frac{\eps}{(x - a_i) + \eps^2}\right)'' f(a_i) \Delta_i
\]
is a Riemann sum for $(P_\eps \ast f)''$. Since $f$ is uniformly continuous, both sets of Riemann sums converge uniformly, so $P_\eps \ast f$ is a limit of such Riemann sums in the norm~\eqref{Eq:Norm-four-terms}.

It remains to note that
\[
\frac{b}{(x - a) + b^2} = \frac{1}{2i} \left( \frac{1}{a + b i - x} - \frac{1}{a - b i - x} \right) \in \mc{D}. \qedhere
\]
\end{proof}

\begin{Prop}
\label{Prop:Continuous}
For a free L{\'e}vy process of the first kind, on each $L^p(\mf{R}, dx)$, $1 \leq p < \infty$ and on $(C_0(\mf{R}), \norm{\cdot}_\infty)$, $\mc{K}_{s,t}$ is strongly continuous in $s, t$.
\end{Prop}

\begin{proof}
Since $\mc{D}$ is dense in $L^p(\mf{R}, \,dx)$ and $\mc{K}_{s,t}$ is a contraction on it, it suffices to prove continuity for $f \in \mc{D}$. Indeed,
\[
\abs{\mc{K}_{s', t'} \left[\frac{1}{z - x} \right] - \mc{K}_{s, t} \left[\frac{1}{z - x} \right]}
= \abs{\frac{F_{s', t'}(z) -  F_{s, t}(z)}{(F_{s', t'}(z) - x)(F_{s, t}(z) - x)}}.
\]
It remains to note that for a fixed $z$, $(F_{s', t'}(z) -  F_{s, t}(z)) \rightarrow 0$ as $s' \rightarrow s$, $t' \rightarrow t$, and that for any $z, w \not \in \mf{R}$,
\begin{equation}
\label{Eq:Quotient-in-Lp}
\norm{\frac{1}{(z - x)(w - x)}}_p \leq \norm{\frac{1}{z - x}}_{2p} \norm{\frac{1}{w - x}}_{2p} \leq C \frac{1}{(\abs{\Im z} \abs{\Im w})^{(1 - p^{-1})/2}}
\end{equation}
for all $p \geq 1$.
\end{proof}

\subsection{Free L{\'e}vy processes of the first kind}
\label{Subsec:FL1}

\begin{Lemma}
\label{Lemma:R-stripped}
Let $\set{\mu_t}$ be a free convolution semigroup, where $\mu = \mu_1$ has mean $\alpha$ and variance $1$. Then
\[
R_\mu(G_{\mu_t}(z)) = \alpha + G_{\nu_t}(z),
\]
where $\nu_t = \mc{J}[\mu_t]$.
\end{Lemma}

\begin{proof}
By definition \eqref{Eq:J-definition},
\[
G_{\nu_t}(z)
= \frac{1}{t} \left( z - \alpha t - \frac{1}{G_{\mu_t}(z)} \right)
= \frac{1}{t} \left( z - \frac{1}{G_{\mu_t}(z)} \right) - \alpha.
\]
One the other hand, by definition of the $R$-transform
\[
G_{\mu_t}^{-1}(z) = \frac{1}{z} + t R_\mu(z),
\]
so
\[
R_{\mu}(z) = \frac{1}{t} \left( G_{\mu_t}^{-1}(z) - \frac{1}{z} \right).
\]
Putting these together, it follows that
\[
\alpha + G_{\nu_t}(z) = R_\mu(G_{\mu_t}(z))
\]
on a domain, and hence, by analytic continuation, on $\mf{C}^+$.
\end{proof}

\begin{Remark}
In \cite{Belinschi-Nica-B_t}, Belinschi and Nica defined a family of transformations
\[
\mf{B}_t[\nu] = \left( \nu^{\boxplus(1 + t)} \right)^{\uplus \frac{1}{1 + t}}.
\]
They showed that these transformations form a semigroup under composition, and $\mf{B}_1 = \mf{B}$ is the Boolean-to-free Bercovici-Pata bijection, defined via
\[
z - F_\nu(z) = z R_{\mf{B}[\nu]}(1/z).
\]
The domain of $\mf{B}$ consists of all probability measures, while its image are all the freely infinitely divisible measures. Moreover, they proved the following evolution equation:
\[
\Phi[\rho \boxplus \gamma_t] = \mf{B}_t[\Phi[\rho]].
\]
We found this equation quite mysterious. We now re-interpret it as follows: a single coefficient stripping, applied to a free convolution semigroup (with finite variance), produces a semicircular evolution started at the free canonical measure of the semigroup.
\end{Remark}

\begin{Prop}
\label{Prop:Evolution}
Let $\rho$ be a probability measure on $\mf{R}$. Then
\[
\mu_t = \Phi_t[\rho \boxplus \gamma_t]
\]
is a free convolution semigroup with mean zero and finite variance, such that $\rho$ is the corresponding free canonical measure. Moreover, each such free convolution semigroup with $\Var[\mu_1] = 1$ arises in this way. In particular, for any such free convolution semigroup,
\[
\mc{J}[\mu_t] = \rho \boxplus \gamma_t.
\]
\end{Prop}

\begin{proof}
We compute
\[
\mu_t = \Phi[\rho \boxplus \gamma_t]^{\uplus t}
= \mf{B}_t[\Phi[\rho]]^{\uplus t}
= \mf{B}_{t-1}[\mf{B}[\Phi[\rho]]]^{\uplus t}
= \mf{B}[\Phi[\rho]]^{\boxplus t},
\]
so $\set{\mu_t}$ form a free convolution semigroup, with $\mu = \mu_1 = \mf{B}[\Phi[\rho]]$. Also,
\[
\mc{J}[\mu_t] = \mc{J}[\Phi_t[\rho \boxplus \gamma_t]] = \rho \boxplus \gamma_t.
\]
If $R_\mu$ is the $R$-transform corresponding to $\set{\mu_t}$, then
\[
G_{\rho \boxplus \gamma_t}(z) = G_{\mc{J}[\mu_t]}(z) = R_\mu(G_{\mu_t}(z))
\]
by the Lemma~\ref{Lemma:R-stripped}. In particular, since $\mu_0 = \delta_0$ and $G_{\mu_0}(z) = \frac{1}{z}$,
\[
G_{\rho}(z) = R_\mu(1/z) = \int_{\mf{R}} \frac{1}{z - x} \,d\rho(x),
\]
so $\rho$ is the free canonical measure for $\set{\mu_t}$. Finally, such a representation holds precisely for any free convolution semigroup with mean zero and $\Var[\mu_1] = 1$.
\end{proof}

\begin{Cor}
For $\set{\mu_t}$ a free convolution semigroup satisfying \eqref{Eq:Normalization}, with free canonical pair $(\alpha, \rho)$,
\[
R_\mu(G_{\mu_s}(z)) = \alpha + G_{\rho \boxplus \gamma_s}(z).
\]
\end{Cor}

\begin{Defn}
For a finite measure $\nu$, we denote by $L_\nu$ the operator
\[
L_\nu[f] = \int_{\mf{R}} \frac{f(x) - f(y)}{x - y} \,d\nu(x) = (I \otimes \nu) (\partial f),
\]
where $\partial$ is the difference quotient operator. Such operators were studied in \cite{Ans-Bochner}, but also in other sources, for example \cite{Liaw-Treil-Rank-one}.
\end{Defn}

\begin{Lemma}
\label{Lemma:A-on-resolvent}
For $\rho$, $\set{\mu_t}$ as in the preceding corollary,
\[
(\alpha \partial_x + \partial_x L_{\rho \boxplus \gamma_s}) \frac{1}{z-x}
= R_\mu(G_{\mu_s}(z)) \frac{1}{(z - x)^2}
\]
and (for each fixed $z \in \mf{C} \setminus \mf{R}$) this is a continuous function of $s$ into $L^p(\mf{R}, \,dx)$.
\end{Lemma}

\begin{proof}
Using the preceding corollary,
\[
(\alpha \partial_x + \partial_x L_{\rho \boxplus \gamma_s}) \frac{1}{z-x} = \alpha \frac{1}{(z - x)^2} + \frac{1}{(z - x)^2} G_{\rho \boxplus \gamma_s}(z)
= R_\mu(G_{\mu_s}(z)) \frac{1}{(z - x)^2}.
\]
Continuity follows from equation~\eqref{Eq:Quotient-in-Lp}.
\end{proof}

\begin{Prop}
\label{Prop:Generator-free}
Let $A_t$ be a generator of a free L{\'e}vy process corresponding to the free convolution semigroup $\set{\mu_t}$ with free canonical pair $(\alpha, \rho)$. Then for $\mc{D}$ from equation~\eqref{Eq:Domain-rational-functions}, $\mc{D} \subset \mc{D}(A_t)$, and on this domain
\begin{equation}
\label{Eq:Generator-free-1st-kind}
A_t f(x) = \alpha \partial_x f(x) + \int_{\mf{R}} \partial_x \frac{f(x) - f(y)}{x - y} \,d(\rho \boxplus \gamma_t)(y),
\end{equation}
which we will abbreviate as
\begin{equation}
\label{Eq:Generator-free-1st-kind-short}
A_t = \alpha \partial_x + \partial_x L_{\rho \boxplus \gamma_t}.
\end{equation}
\end{Prop}

\begin{proof}
For any free convolution semigroup $\set{\mu_t}$, the following evolution equation holds:
\begin{equation}
\label{Eq:Evolution}
\partial_t G_{\mu_t}(z) = - R_\mu(G_{\mu_t}(z)) \ G_{\mu_t}'(z),
\end{equation}
see equation~(3.18) in \cite{VDN}. Also according to Theorem~3.1 of \cite{BiaProcesses}, the transition operators of a free L{\'e}vy process (of the first kind) have the property that
\[
\mc{K}_{s,t} \left[ \frac{1}{z - x} \right] = \frac{1}{F_{s,t} - x},
\]
where
\[
G_{\mu_t}(z) = G_{\mu_s}(F_{s,t}(z)).
\]
This implies that
\[
\partial_t F_{s,t}(z) = \frac{\partial_t G_{\mu_t}(z)}{G_{\mu_s}'(F_{s,t}(z))}
= - \frac{G_{\mu_t}'(z) R(G_{\mu_t}(z))}{G_{\mu_s}'(F_{s,t}(z))}.
\]
So using Lemma~\ref{Lemma:A-on-resolvent}, we compute
\[
\begin{split}
& \abs{\frac{1}{h} \left(\mc{K}_{s, s+h} \frac{1}{z - x} - \frac{1}{z - x} \right) - \left(\alpha \partial_x + \partial_x L_{\rho \boxplus \gamma_s} \right) \frac{1}{z - x}} \\
&\qquad = \abs{ \frac{1}{h} \left( \frac{1}{F_{s, s+h}(z) - x} - \frac{1}{z - x} \right) - \frac{1}{(z-x)^2} R_\mu(G_{\mu_s}(z))} \\
&\qquad = \abs{- \frac{1}{h} \frac{F_{s, s+h}(z) - F_{s,s}(z)}{(F_{s,s+h}(z) - x) (z-x)} - \frac{1}{(z-x)^2} R_\mu(G_{\mu_s}(z))} \\
&\qquad = \abs{\left( \frac{1}{h} \int_0^h \frac{G_{\mu_{s+u}}'(z)}{G_{\mu_s}'(F_{s,s+u}(z))} \frac{R_\mu(G_{\mu_{s+u}}(z))}{(F_{s,s+h}(z) - x)} \,du - \frac{R_\mu(G_{\mu_s}(z))}{(z-x)} \right) \frac{1}{z-x}} \\
&\qquad = \abs{\left( \frac{1}{h} \int_0^h \frac{G_{\mu_{s+u}}'(z)}{G_{\mu_s}'(F_{s,s+u}(z))} \frac{R_\mu(G_{\mu_{s+u}}(z))}{R_\mu(G_{\mu_s}(z))} \,du \right) \frac{1}{(F_{s,s+h}(z) - x)} - \frac{1}{(z-x)}} \abs{\frac{R_\mu(G_{\mu_s}(z))}{z-x}}.
\end{split}
\]
Now using equation~\eqref{Eq:Quotient-in-Lp} and $\frac{G_{\mu_{s+u}}'(z)}{G_{\mu_s}'(F_{s,s+u}(z))} \frac{R_\mu(G_{\mu_{s+u}}(z))}{R_\mu(G_{\mu_s}(z))} \rightarrow 1$ as $u \rightarrow 0$, the difference above converges to zero in $L^p(\mf{R}, \,dx)$, $1 \leq p \leq \infty$. The result follows.
\end{proof}

The appearance of the semicircular distributions in the generator formula above is explained by the following theorem.

\begin{Thm}
Let $\set{X_t : t \geq 0}$ be a centered free L{\'e}vy process of the first kind with finite variance, normalized so that $\Var[X_1] = 1$. Let $\rho$ be the canonical measure for the corresponding free convolution semigroup $\set{\mu_t}$. Finally, let $\set{Y_t : t \geq 0}$ be the free Brownian motion started at $Y_0$ with distribution $\rho$. Then the transition operators of the processes $\set{X_t}$ and $\set{Y_t}$ coincide.
\end{Thm}

\begin{proof}
It suffices to check the equality of transition operators on $\mc{D}$, in other words we need to verify the equality of analytic functions $F_{s,t}$. We check that indeed,
\[
\begin{split}
G_{\mu_s}^{-1} \circ G_{\mu_t}(z)
& = \Bigl(G_{\mu_t}^{-1} - (t-s) R_\mu \Bigr) \circ G_{\mu_t}(z)
= z - (t-s) R_\mu \circ G_{\mu_t}(z) \\
& = z - (t-s) G_{\rho \boxplus \gamma_t}(z)
= \Bigl(G_{\rho \boxplus \gamma_t}^{-1} - (t-s) z \Bigr) \circ G_{\rho \boxplus \gamma_t}(z)
= G_{\rho \boxplus \gamma_s}^{-1} \circ G_{\rho \boxplus \gamma_t}(z).
\end{split}
\]
\end{proof}

\begin{Remark}
For readers familiar with the properties of the subordination distribution, we provide an alternative proof, see \cite{Nica-Subordination} for the results used. We compute
\[
\begin{split}
\mu_{t-s} \subord \mu_s
& = (\mu \subord \mu_s)^{\boxplus(t-s)}
= (\mu_s^{\boxplus (1/s)} \subord \mu_s)^{\boxplus (t-s)}
= \mf{B}[\mu_s]^{\boxplus (1/s)(t-s)}
= \mf{B}\left[\left(\mu^{\boxplus s}\right)^{\uplus (1/s)}\right]^{\boxplus (t-s)} \\
& = \left( \mf{B} \circ \mf{B}_{s-1} [\mu] \right)^{\boxplus (t-s)}
= \left( \mf{B}_{s} \circ \mf{B}[\Phi[\rho]] \right)^{\boxplus (t-s)}
= \left(\mf{B} \circ \mf{B}_s[\Phi[\rho]] \right)^{\boxplus (t-s)} \\
& = \mf{B}[\Phi[\rho \boxplus \gamma_s]]^{\boxplus (t-s)}
= \Bigl(\gamma \subord (\rho \boxplus \gamma_s) \Bigr)^{\boxplus (t-s)}
= \gamma_{t-s} \subord (\rho \boxplus \gamma_s).
\end{split}
\]
Note also that the preceding theorem is false for a process with non-zero mean; indeed, the generator of a free Brownian motion with drift is not \eqref{Eq:Generator-free-1st-kind-short} but rather $\alpha \partial_x + \partial_x L_{\rho \boxplus \gamma_t \boxplus \delta_{\alpha t}}$.
\end{Remark}

\begin{Prop}
\label{Prop:Lp-estimates}
Let $\nu$ be a finite measure. Then
\[
\norm{\partial_x L_\nu f}_\infty \leq \nu(\mf{R}) \norm{f''}_\infty.
\]
and for $1 \leq p < \infty$,
\[
\norm{\partial_x L_\nu f}_1 \leq C_p \nu(\mf{R}) (\norm{f''}_\infty + \norm{f'}_p).
\]
It follows that $\partial_x L_\nu$ is a bounded operator $W^\infty \rightarrow C_0(\mf{R})$ and $W^p \rightarrow L^p(\mf{R}, \,dx)$.
\end{Prop}

\begin{proof}
By Taylor's theorem,
\[
\begin{split}
\abs{\partial_x \frac{f(x) - f(y)}{x - y}}
& = \abs{\frac{f(y) - f(x) - (y-x) f'(x)}{(y-x)^2}} \\
& = \abs{\frac{1}{(y - x)^2} \int_x^y (y - u) f''(u) \,du}
\leq \sup_{x \leq u \leq y} \abs{\frac{y - u}{y - x} f''(u)} \leq \norm{f''}_\infty.
\end{split}
\]
So
\begin{equation}
\label{Eq:Norm-second-derivative}
\norm{\partial_x L_\nu f}_\infty \leq \nu(\mf{R}) \norm{f''}_\infty.
\end{equation}
Since $\partial_x L_\nu(\mc{D}) \subset C_0(\mf{R})$, $\mc{D}$ is dense in $W^\infty$, and $C_0(\mf{R})$ is closed, it follows that in fact
\[
\partial_x L_\nu(W^\infty) \subset C_0(\mf{R}).
\]

Next, \eqref{Eq:Norm-second-derivative} implies that for $p > 1$ and $q$ the dual exponent,
\[
\begin{split}
\norm{\partial_x L_\nu f}_1
& = \int_{\mf{R}} \int_{[y - a, y + a]} \abs{ \partial_x \frac{f(x) - f(y)}{x - y}} \,dx \,d\nu(y) \\
&\quad + \int_{\mf{R}} \int_{[y - a, y + a]^c} \abs{ \frac{f'(x)}{x - y} - \frac{f(x) - f(y)}{(x - y)^2}} \,dx \,d\nu(y) \\
& \leq \int_{\mf{R}} 2 a \norm{f''}_\infty \,d\nu(y) + \int_{\mf{R}} \norm{f'}_p \left( \int_{[y - a, y + a]^c} \frac{1}{(x - y)^q} \,dx \right)^{1/q} \,d\nu(y) \\
&\quad + \int_{\mf{R}} \int_{[y - a, y + a]^c} \frac{1}{(x - y)^2} \left( \int_y^x \,du \right)^{1/q} \left( \abs{f'(u)}^p \,du \right)^{1/p} \,dx \,d\nu(y) \\
& \leq 2 a \nu(\mf{R}) \norm{f''}_\infty + \frac{(2 (p-1))^{1/q}}{a^{1/p}} \nu(\mf{R}) \norm{f'}_p + \frac{2 p}{a^{1/p}} \nu(\mf{R}) \norm{f'}_p
\end{split}
\]
So
\[
\norm{\partial_x L_\nu f}_1 \leq C_p \nu(\mf{R}) \left( \norm{f''}_{\infty} + \norm{f'}_{p} \right).
\]
A similar argument works for $p=1$. The final result follows by interpolation.
\end{proof}

\begin{Thm}
\label{Thm:Generator-free}
Let $\rho$ be a probability measure, $\alpha \in \mf{R}$, $\set{\mu_t}$ the free convolution semigroup with the free canonical pair $(\alpha, \rho)$, $\set{X(t)}$ the corresponding free L{\'e}vy process, and $\set{A_t}$ its generators. Then on $C_0(\mf{R})$, $(A_t, W^\infty)$ equals
\begin{equation*}
A_t = \alpha \partial_x + \partial_x L_{\rho \boxplus \gamma_t}.
\end{equation*}
and on $L^p(\mf{R}, \,dx)$, $(A_t, W^p)$ is given by the same formula.
\end{Thm}

\begin{proof}
Use the estimates in Proposition~\ref{Prop:Lp-estimates}, and Lemma~\ref{Lemma:W^p-density}, and apply Proposition~\ref{Prop:Basic} with $\mc{A} = C_0(\mf{R})$, $\mc{B} = W^\infty$, respectively, $\mc{A} = L^p(\mf{R}, \,dx)$, $\mc{B} = W^p$, to show that these sets are in the domain of the generators. Since the same estimate shows that $\alpha \partial_x + \partial_x L_{\rho \boxplus \gamma_t}$ continuously extends to $\mc{B}$, the result follows.
\end{proof}

\begin{Example}
A free Meixner distribution $\mu_{b,c}^{\boxplus t}$ is the probability measure with Jacobi parameters
\[
\alpha_0 = 0, \alpha_n = b, \beta_1 = t, \beta_n = t + c.
\]
For $c \geq 0$, these distributions form a free convolution semigroup with respect to the parameter $t$. Clearly the corresponding $\nu_t = \mc{J}[\mu_{b,c}^{\boxplus t}] = \delta_{b} \boxplus \gamma_{t + c}$ are the semicircular distributions with mean $b$ and variance $(t+c)$; thus we recover a weaker version of the result of \cite{Bryc-Markov-Meixner}. On the other hand, for $\mu = \mu_{b,c}$ we also have
\[
R_\mu(z) = z \left(1 + b R_\mu(z) + c (R_\mu(z))^2 \right),
\]
which implies that
\[
R_{\mu}(z) = \int_{\mf{R}} \frac{z}{1 - z x} \,d\rho(x)
\]
for $\rho = \delta_{b} \boxplus \gamma_c$ semicircular with mean $b$ and variance $c$. So the free canonical measure in this case is also semicircular. The reason for this coincidence is that, as pointed out in Proposition~\ref{Prop:Evolution},
\[
\nu_t = \rho \boxplus \gamma_t = \left( \delta_{b} \boxplus \gamma_c \right) \boxplus \gamma_t = \delta_{b} \boxplus \gamma_{t + c}.
\]

In the particular case $b = c = 0$, we have $\mu_t = \nu_t = \gamma_t$ (and $\rho = \delta_0$). The corresponding process is the free Brownian motion, whose generator $\partial_x L_{\gamma_t}$ was found at the end of Section~4 in \cite{BiaSpeBrownian}, see also Example~4.9 in \cite{BKSQGauss}.
\end{Example}

\begin{Example}[Generator of the Cauchy process]
The mean of the Cauchy distribution is undefined. Nevertheless, we can still compute the generator of the free Cauchy process, because the Cauchy distributions form both a free and a usual convolution semigroup. Indeed, the Fourier transform of the standard Cauchy distribution is $e^{-\abs{x}}$, so the generator of the corresponding process is $- \abs{-i \partial_x} = - \abs{\partial_x}$. Note that $\abs{x} = \sgn(x) x$ and
\[
\mc{F}(H f)(x) = - i \sgn(x) \mc{F}(f)(x),
\]
where $\mc{F}$ is the Fourier transform and $H$ is the Hilbert transform. We conclude that the generator is
\[
A f(x) = -i \partial_x (H f)(x).
\]
This is consistent with the relation $R_\mu(z) = -i$ and
\[
\partial_t G_{\mu_t}(z) = - i G_{\mu_t}'(z).
\]
\end{Example}

\begin{Remark}
The generator of a classical process is a pseudo-differential operator. The generator of the free process can be given a similar interpretation. Indeed, note first that
\[
L_{\delta_0} f(x) = \frac{f(x) - f(0)}{x}.
\]
This operator is the crucial object in \cite{AnsBoolean}; note also that $\partial_x L_{\delta_0}$ is the generator, but only at time zero, of the free Brownian motion. Suppose now that all the moments $m_n(\nu)$ of $\nu$ are finite, and let
\[
M_\nu(z) = \sum_{n=0}^\infty m_n(\nu) z^n
\]
be its moment generating function. Then, at least for polynomial $f$,
\[
L_\nu f = M_\nu(L_{\delta_0}) L_{\delta_0} f.
\]
Indeed, for $f(x) = x^n$,
\[
M_\nu(L_{\delta_0}) L_{\delta_0} f(x)
= \sum_{k=0}^n m_k(\nu) L_{\delta_0}^{k+1} x^n
= \sum_{k=0}^{n-1} m_k(\nu) x^{n-k-1}
= \int_{\mf{R}} \frac{x^n - y^n}{x - y} \,d\nu(y)
= L_\nu f(x).
\]
By writing $L_\nu = G_\nu(L_{\delta_0}^{-1})$, we can interpret it as a pseudo-differential-type operator even if the moments of $\nu$ are not finite.
\end{Remark}

\begin{Remark}
Suppose $\mu$ is compactly supported. In particular, we can expand
\[
R_{\mu}(z) = \sum_{n=1}^\infty r_n z^{k-1},
\]
where $\set{r_n}$ are the free cumulants of $\mu$. Then it follows from Proposition~\ref{Prop:Generator-free}, Lemma~\ref{Lemma:R-stripped}, and the preceding remark that the generator of the corresponding process is
\begin{equation}
\label{Eq:Generator-sum}
\begin{split}
A_s & = \alpha \partial_x + \partial_x L_{\rho \boxplus \gamma_s}
= \partial_x \Bigl(\alpha + G_{\rho \boxplus \gamma_s}(L_{\delta_0}^{-1}) \Bigr) \\
& = \partial_x R_\mu(G_{\mu_s}(L_{\delta_0}^{-1}))
= \partial_x R_\mu(L_{\mu_s})
= \sum_{n=1}^\infty r_n \partial_x L_{\mu_s}^{n-1}.
\end{split}
\end{equation}
On the other hand, by Lemmas~2 and 3 of \cite{AnsIto}, in this case the higher diagonal measures
\[
\Delta_n(t) = \int_0^t (d X(s))^n
\]
are defined, and $\Exp{\Delta_n(t)} = r_n t$. Moreover, by Corollary~12 of the same paper, for polynomial $f$
\[
f(X(t)) = \sum_{n=1}^\infty \int_0^t (I \otimes \mf{E} \otimes \ldots \otimes \mf{E} \otimes I) \left[(\partial^n f)(X(s), \ldots, X(s)) \right] \sharp \,d\Delta_n(s),
\]
where $\partial^n$ is defined recursively by
\[
\partial^n = (I \otimes \ldots \otimes I \otimes \partial) \partial^{n-1}
\]
(this notation differs by a factor of $n!$ from \cite{AnsIto}), and
\[
\int_0^t (A(x) \otimes B(s)) \sharp \,dX(s) = \int _0^t A(x) \,dX(s) \, B(s).
\]
In other words,
\[
f(X(t)) = \sum_{n=1}^\infty \int_0^t (\partial L_{\mu_s}^{n-1} f)(X(s)) \sharp \,d\Delta_n(s).
\]

It follows that for the generator \eqref{Eq:Generator-sum}, the martingale from Lemma~\ref{Lemma:Generator-martingale} can be written explicitly as
\[
\begin{split}
& f(X(t)) - \int_0^t (A_s f)(X(s)) \,ds \\
&\qquad = \sum_{n=1}^\infty \int_0^t (\partial L_{\mu_s}^{n-1} f)(X(s)) \sharp \,d\Delta_n(s) - \sum_{n=1}^\infty \int_0^t r_n (\partial_x L_{\mu_s}^{n-1} f) (X(s)) \,ds \\
&\qquad = \sum_{n=1}^\infty \int_0^t (\partial L_{\mu_s}^{n-1} f)(X(s)) \,\sharp \,d(\Delta_n(s) - r_n s).
\end{split}
\]
It would be interesting to find such a representation for more general $\mu$.
\end{Remark}

\begin{Remark}
A matricial interpolation between classical and free L{\'e}vy processes was constructed in \cite{Cabanal-Matrix}, where the generator of such a matricial process is also computed.
\end{Remark}

\subsection{Free L{\'e}vy processes of the second kind}
\label{Subsec:FL2}

The semigroup $\set{\mc{K}_t}$ of transition operators corresponding to a free L{\'e}vy process of the second kind is characterized by
\[
\mc{K}_t \left[ \frac{1}{z - x} \right] = \frac{1}{F_{\mu_t}(z) - x},
\]
and $\set{F_{\mu_t}}$ form a semigroup with respect to composition.

The transition operators for a monotone L{\'e}vy process corresponding to the family $\set{\mu_t}$ are exactly the same, and in fact $\set{\mu_t}$ form a monotone convolution semigroup, see Corollary~5.3 of \cite{Franz-Monotone-Boolean}. According to Theorem~5.1 of \cite{Franz-Muraki-Markov-monotone}, at least for the compactly supported case, on bounded continuous functions
\[
\mc{K}_t f (x) = \int_{\mf{R}} f(y) \,d(\delta_x \triangleright \mu_t)(y).
\]
Since $\set{\mc{K}_t}$ form a semigroup, we only need to compute the generator $A$ at zero. By Proposition~5.1 of \cite{Franz-Muraki-Markov-monotone}, the generator is
\begin{equation}
\label{Eq:Generator-2nd-kind}
A f (x) = \alpha \partial_x f(x) + \int_{\mf{R}} \partial_x \frac{f(x) - f(y)}{x - y} \,d\rho(y),
\end{equation}
where
\[
- \left. \frac{\partial}{\partial t} \right|_{t=0} F_{\mu_t}(z) = z^2 \left. \frac{\partial}{\partial t} \right|_{t=0} G_{\mu_t}(z) = \alpha + \int_{\mf{R}} \frac{1}{z - x} \,d\rho(x),
\]
so that $(\alpha, \rho)$ is the monotone canonical pair (note our choice of signs is the opposite of \cite{Franz-Muraki-Markov-monotone}). As pointed out in Theorem~4.5 of \cite{BiaProcesses}, only certain $\rho$ correspond to processes with free increments in this way (note that unlike Biane, we have assumed $\mu_0 = \delta_0$). We repeat Biane's question (Section~4.7): it would be interesting to have a more direct description of which $\rho$ do so appear. In particular, according to J.C.~Wang, a centered FL2 process cannot have finite variance, and as of this writing, no non-trivial examples of FL2 processes are known.

\begin{Remark}
Let $\set{X(t)}$ be a process whose increments are stationary and independent in a certain sense, $\set{\mc{K}_{s,t}}$ the corresponding transition operators, and $\set{\mu_t, \star}$ be the corresponding convolution semigroup. Since $\mu_0 = \delta_0$, we observe that
\[
\mu_t[f] = \Exp{f(X(t))} = \Exp{(\mc{K}_{0,t} f)(X(0))} = \mu_0[\mc{K}_{0,t} f] = (\mc{K}_{0,t} f)(0).
\]
Therefore the corresponding cumulant functional \eqref{Eq:Cumulant-definition} is
\[
C_\mu[f] = \left. \frac{\partial}{\partial t} \right|_{t=0} \mu_t[f] = (A_0 f)(0).
\]
We note that indeed, if we set $t=0$ and $x=0$ in formulas \eqref{Eq:Generator-classical}, \eqref{Eq:Generator-free-1st-kind} and \eqref{Eq:Generator-2nd-kind}, in all three cases we get formula \eqref{Eq:Cumulant}. Note that in these three cases, $\rho$ is interpreted as the classical, free, and monotone canonical measure, respectively. In particular, in all three cases, the cumulant functionals are defined at least on the domain of the corresponding generators. On the other hand, for $t > 0$
\[
L_t[f] = \partial_t \mu_t[f] = (\mc{K}_{0,t} A_t f)(0) = \mu_t[A_t f]
\]
will depend on the type of semigroup considered.
\end{Remark}

\subsection{Semigroups for generators of the free L{\'e}vy processes of the first kind}
\label{Subsec:Semigroups}

We noted before that for a L{\'e}vy process of the first kind, $(A_t, \mc{D}(A_t))$ is closable. We now show explicitly that its closure generates a contraction semigroup.

Let $\alpha \in \mf{R}$, and $\rho$ be a probability measure. For each $s \geq 0$, denote $\set{F^{(s)}_t : t \geq 0}$ the solution of
\[
- \partial_t F^{(s)}_t(z) = \alpha + G_{\rho \boxplus \gamma_s}(F^{(s)}_t(z))
\]
which, by Theorem~4.5 of \cite{BiaProcesses}, exists and moreover satisfies
\[
F^{(s)}_t(z) = F_{\tau^{(s)}_t}(z)
\]
for a probability measure $\tau^{(s)}_t$. In fact $F^{(s)}_{t_1} \circ F^{(s)}_{t_2} = F^{(s)}_{t_1 + t_2}$, and the corresponding measures form a monotone convolution semigroup. Denote now
\[
\mc{L}^{(s)}_t f(x) = \int_{\mf{R}} f(y) d (\delta_x \triangleright \tau^{(s)}_t)(y).
\]
At least in the compactly supported case, as noted above, these are transition operators for the corresponding process with monotone independent increments, but will not use this property directly. Instead, we note the following.

\begin{Thm}
Let $\set{\mc{L}^{(s)}_t}$ be as above, and $\set{\mc{K}_{s,t}}$, $\set{A_t}$ be as in Theorem~\ref{Thm:Generator-free}.
\begin{enumerate}
\item
For each $s$, the operators $\set{\mc{L}^{(s)}_t}$ form a strongly continuous semigroup of contractions on $C_0(\mf{R})$ and $L^p(\mf{R}, \,dx)$.
\item
The generator $B_s$ of this semigroup is a closed operator for which $\mc{D}$ is a core.
\item
For $f \in \mc{D}$, $B_s f = A_s f = (\alpha \partial_x + \partial_x L_{\rho \boxplus \gamma_s}) f$.
\item
On $C_0(\mf{R})$ and $L^p(\mf{R}, \,dx)$.
\begin{equation}
\label{Eq:Limit}
\lim_{n \rightarrow \infty} \left( \mc{K}_{s, s + t/n} \right)^n = \mc{L}^{(s)}_t
\end{equation}
strongly.
\item
$B_s = \overline{A_s}$, and $\mc{D}$ is a core for $A_s$.
\end{enumerate}
\end{Thm}

\begin{proof}
The proofs are very similar to the results for $\mc{K}_{s,t}$, and are mostly omitted. For part (a), see Propositions~\ref{Prop:Contraction} and \ref{Prop:Continuous}. For the semigroup property, we compute
\[
\begin{split}
\mc{L}^{(s)}_{t_1} \mc{L}^{(s)}_{t_2} f(x)
& = \int_{\mf{R}} \left( \int_{\mf{R}} f(z) d(\delta_y \triangleright \tau^{(s)}_{t_2})(z) \right) d(\delta_x \triangleright \tau^{(s)}_{t_1})(y) \\
& = \int_{\mf{R}} f(z) d((\delta_x \triangleright \tau^{(s)}_{t_1}) \triangleright \tau^{(s)}_{t_2})(z) \\
& = \int_{\mf{R}} f(z) d(\delta_x \triangleright \tau^{(s)}_{t_1 + t_2})(z)
= \mc{L}^{(s)}_{t_1 + t_2} f(x),
\end{split}
\]
since $\triangleright$ is associative and distributive in the first variable.

The generator of a strongly continuous semigroup is closed. $\mc{D}$ is dense and invariant under all $\mc{L}^{(s)}_t$, so by Theorem X.49 of \cite{ReeSim2} it is a core for $B_s$. The proof of part (c) is similar to Proposition~\ref{Prop:Generator-free}. Part (d) follows from Chernoff Product Formula, Theorem~III.5.2 of \cite{Engel-Nagel-book}, applied to $(A_s, \mc{D}) = (B_s, \mc{D})$ (the density of the range condition is satisfied since $B_s$ generates a contraction semigroup).

Finally, for part (e), we apply Chernoff Product Formula to $(A_s, \mc{D}(A_s))$. The density of the range condition holds since it was already satisfied for $B_s$ on $\mc{D}$. The theorem implies that $\overline{A_s}$ generates precisely the same semigroup~\eqref{Eq:Limit}. Therefore $\overline{A_s} = B_s$.
\end{proof}

\begin{Remark}
A alternative general approach to the non-autonomous Cauchy problem is via evolution semigroups, see Section VI.9(b) of \cite{Engel-Nagel-book}. We briefly describe this approach in our case. Let $f \in C_0(\mf{R}^2)$, and denote $f_t(x) = f(t,x)$. Then operators
\[
(T_t f) (s,x) = (\mc{K}_{s,s+t} f_{s+t})(x)
\]
form a contraction semigroup with respect to $\norm{\cdot}_\infty$. Its generator is a closed, dissipative operator. At least formally, it is related to the generators of the family $\set{\mc{K}_{s,t}}$ by
\[
(\mb{A} f)(t,x) = A_t f_t(x) + \partial_t f_t(x).
\]
We also note that $\mb{A} f = 0$ if $f$ is fixed by $T_t$, in other words if
\[
\mc{K}_{s,s+t} f_{s+t} = f_s.
\]
This is precisely the condition for $f_t(X_t)$ to be a martingale.
\end{Remark}

\subsection{Further remarks on the properties of $L_\nu$ and $\partial_x L_\nu$}
\label{Subsec:Further}

\begin{Prop}
Let $\mu$ be a probability measure with mean $\alpha$ and variance $\beta$. Denote $\nu = \mc{J}[\mu]$. Then $L_\mu$ is a multiple of a unitary operator
\[
\set{f \in L^2(\mu): \mu[f] = 0} \rightarrow L^2(\nu),
\]
with inverse $x - \alpha - \beta L_\nu$.
\end{Prop}

Note that if polynomials are dense in $L^2(\mu)$, this result follows from the proof of Proposition~10 in \cite{Ans-Bochner}, and the statement about the inverse from that proposition and Corollary~11.

\begin{proof}
First we show that
\[
\widetilde{\mc{D}} = \Span{ \frac{1}{z - x} : z \in \mf{C} \setminus \mf{R}}
\]
is dense in $L^2(\mu)$. Indeed, if $f \in \widetilde{\mc{D}}^\perp$, then
\[
\ip{f}{\frac{1}{\bar{z} - x}}_\mu = \int_{\mf{R}} \frac{f(x)}{z - x} \,d\mu(x) = G_{f \mu}(z) = 0
\]
for all $z \in \mf{C} \setminus \mf{R}$. By Stieltjes inversion, it follows that $f  = 0$ $\mu$-a.e.

By density of $\widetilde{\mc{D}}$, it suffices to show that for resolvents,
\[
\begin{split}
\ip{\frac{1}{z - x} - G_\mu(z)}{\frac{1}{w - x} - G_\mu(w)}_\mu
& = \mu \left[ \frac{1}{z - x} \frac{1}{\overline{w} - x} \right] - G_\mu(z) G_\mu(\overline{w}) \\
& = \frac{1}{z - \overline{w}} (G_\mu(\overline{w}) - G_\mu(z)) - G_\mu(z) G_\mu(\overline{w}) \\
& = \frac{G_\mu(z) G_\mu(\overline{w})}{\overline{w} - z} \left( (z - F_\mu(z)) - (\overline{w} - F_\mu(\overline{w})) \right) \\
& = \frac{G_\mu(z) G_\mu(\overline{w})}{\overline{w} - z} (\beta G_\nu(z) - \beta G_\nu(\overline{w})) \\
& = \beta \ip{ \frac{1}{z - x} G_\mu(z)}{\frac{1}{w - x} G_\mu(w)}_\nu \\
& = \beta \ip{L_\mu \left[ \frac{1}{z - x} \right]}{L_\mu \left[ \frac{1}{w - x} \right]}_\nu.
\end{split}
\]
independently of $\alpha$. Also
\[
(x - \alpha - \beta L_\nu) \circ L_\mu \left[ \frac{1}{z - x} \right]
= \left( -(z - x) + z - \alpha - \beta G_\nu(z) \right) \frac{1}{z - x} G_\mu(z)
= - G_\mu(z) + \frac{1}{z - x}.
\]
and the proof of the opposite equality is similar.
\end{proof}

We noted in Proposition~\ref{Prop:Basic} combined with Theorem~\ref{Thm:Generator-free} that a generator of the transition operator family for a free L{\'e}vy process is dissipative, and that on $W^p$, it coincides with $\alpha \partial_x + \partial_x L_\nu$. We now give a more explicit proof of this result for $p=2$, which may be of interest in itself.

\begin{Prop}
Let $\nu$ be a finite measure, and denote $A = \alpha \partial_x + \partial_x L_\nu$ and
\[
D(f, \bar{g}) = \int_{\mf{R}} (\partial f)(x,y) (\partial g)(x,y) \,d\nu(y).
\]
Then $D$ is well-defined for $f, g \in W^{1, \infty} \cap W^{1, 2}$, while for $f, g \in W^2$
\[
D(f, \bar{g}) = A(f g) - f A(g) - A(f) g
\]
for all $\alpha$, so that $D$ is the carr{\'e} du champ operator corresponding to $A$. It follows that for such $f$,
\[
\Re \ip{A f}{f} \leq 0,
\]
so $\alpha \partial_x + \partial_x L_\nu$ on $W^2$ is $L^2$-dissipative.
\end{Prop}

\begin{proof}
We compute
\[
\begin{split}
\int_{\mf{R}} D(f, f) \,dx
& = \iint_{\mf{R}^2} \abs{\frac{f(x) - f(y)}{x - y}}^2 \,d\nu(y) \,dx \\
& = \int_{\mf{R}} \int_{[y - a, y + a]} \abs{\frac{f(x) - f(y)}{x - y}}^2  \,dx \,d\nu(y) + \int_{\mf{R}} \int_{[y - a, y + a]^c} \abs{\frac{f(x) - f(y)}{x - y}}^2 \,dx \,d\nu(y) \\
& \leq \int_{\mf{R}} 2 a \norm{f'}_\infty^2 \,d\nu(y) + \int_{\mf{R}} \int_{[y - a, y + a]^c} \frac{1}{(x - y)^2} \abs{\int_y^x f'(u) \,du}^2 \,dx \,d\nu(y) \\
& \leq 2 a \norm{f'}_\infty^2 \nu(\mf{R}) + \norm{f'}_1^2 \int_{\mf{R}} \frac{2}{a} \,d\nu(y)
= 2 \nu(\mf{R}) \left( a \norm{f'}_\infty^2 + \frac{1}{a} \norm{f'}_1^2 \right).
\end{split}
\]
For $a = 1$, we get
\[
\norm{D(f, f)}_1 \leq 2 \nu(\mf{R}) \left( \norm{f}_{1, \infty} + \norm{f}_{1,1} \right)^2.
\]
By polarization, $D(f,g)$ is well defined for $f, g \in W^{1, \infty} \cap W^{1,1}$.

Next, we compute
\[
\begin{split}
\left(A(f g) - f A(g) - A(f) g \right) (x)& = \int_{\mf{R}} \left[ \frac{f'(x) g(x) - f(x) g'(x)}{x - y} - \frac{f(x) g(x) - f(y) g(y)}{(x - y)^2} \right] \,d\nu(y) \\
&\quad - \int_{\mf{R}} \left[ \frac{f'(x) g(x)}{x - y} - \frac{f(x) - f(y)}{(x - y)^2} g(x) \right] \,d\nu(y) \\
&\quad - \int_{\mf{R}} \left[ \frac{f(x) g'(x)}{x - y} - f(x) \frac{g(x) - g(y)}{(x - y)^2} \right] \,d\nu(y) \\
&\quad + \alpha ((fg)' - f g' - f' g)(x) \\
& = \int_{\mf{R}} \frac{f(x) - f(y)}{x - y} \frac{g(x) - g(y)}{x - y} \,d\nu(y).
\end{split}
\]
Moreover,
\[
\begin{split}
\abs{L_\nu(f)(x)}
& = \int_{[x - a, x + a]} \frac{1}{\abs{x - y}} \abs{\int_y^x f'(u) \,du} \,d\nu(y) +  \int_{[x - a, x + a]^c} \frac{1}{\abs{x - y}} \abs{\int_y^x f'(u) \,du} \,d\nu(y) \\
& \leq \norm{f'}_\infty \nu([x - a, x + a]) \\
&\quad + \int_{[x - a, x + a]^c} \frac{1}{\abs{x - y}} \left( \int_y^x \,du \right)^{1/2} \abs{ \int_y^x \abs{f'(u)}^2 \,du}^{1/2} \,d\nu(y) \\
& \leq \norm{f'}_\infty \nu([x - a, x + a]) + \frac{1}{\sqrt{a}} \nu(\mf{R}) \norm{f'}_2.
\end{split}
\]
So for $a = 1/\eps^2$,
\[
\limsup_{x \rightarrow \infty} \abs{L_\nu(f)(x)} \leq \eps \nu(\mf{R}) \norm{f'}_2
\]
as long as $\norm{f'}_\infty < \infty$. Thus for $f \in W^{1, \infty} \cap W^{1, 2}$,
\begin{equation}
\label{Eq:Limit-zero}
\limsup_{x \rightarrow \infty} \abs{L_\nu(f)(x)} = 0.
\end{equation}
Therefore for $f \in W^2$,
\begin{equation}
\label{Eq:Dissipative}
\begin{split}
2 \Re \ip{A f}{f} = \ip{A f}{f} + \ip{f}{A f}
& = \int_{-\infty}^{\infty} \left[A(\abs{f}^2)(x) - D(f, \bar{f})(x) \right] \,dx \\
& = \int_{-\infty}^{\infty} \left[ \partial_x \Bigl( L_\nu(\abs{f}^2)(x) + \alpha f(x) \Bigr)- D(f, \bar{f})(x) \right] \,dx \\
& = - \iint_{\mf{R}^2} \abs{\frac{f(x) - f(y)}{x - y}}^2 \,d\nu(y) \,dx < 0.
\end{split}
\end{equation}
\end{proof}

Note that the Dirichlet form from \cite{Biane-Log-Sobolev} is $\mc{D}(f, g) = \int_{\mf{R}} D(f, g)(x) \,d\nu(x)$, which is not the same as the right-hand-side in equation~\eqref{Eq:Dissipative}.

\begin{Prop}
$\alpha \partial_x + \partial_x L_\nu$ is $C_0$-dissipative on $W^\infty$.
\end{Prop}

\begin{proof}
For $f \in C_0(\mf{R})$ such that $\abs{f(x_0)} = \max_{x \in \mf{R}} \abs{f(x)}$, a normalized tangent functional is $\phi_f(g) = \overline{f(x_0)} \delta_{x_0}$. Then for $f \in W^\infty$,
\[
\begin{split}
\Re \phi_f[\alpha \partial_x f + \partial_x L_\nu f]
& = \Re \delta_{x_0} \left[ \alpha \overline{f(x_0)} f'(x) + \overline{f(x_0)} \int_{\mf{R}} \left( \frac{f'(x)}{x-y} - \frac{f(x) - f(y)}{(x-y)^2} \right) \,d\nu(y) \right] \\
& = - \Re \left( \overline{f(x_0)} \int_{\mf{R}} \frac{f(x_0) - f(y)}{(x_0 - y)^2} \right) \,d\nu(y) \\
& = - \int_{\mf{R}} \frac{\abs{f(x_0)}^2 - \Re (\overline{f(x_0)} f(y))}{(x_0 - y)^2} \,d\nu(y)
\leq 0.
\end{split}
\]
since
\[
2 \Re(\overline{f(x_0)} f'(x_0)) = (\overline{f} f)'(x_0) = (\abs{f}^2)'(x_0) = 0.
\]
\end{proof}

\section{The $q$-Brownian motion}
\label{Section:q}

Let $q \in (-1,1)$. The $q$-Brownian motion $\set{X(t) : t \geq 0}$ is a non-commutative stochastic process constructed in \cite{BozSpeBM1}. The distribution of each $X(t)$ is a very classical $q$-Gaussian distribution
\begin{equation}
\label{Eq:q-Gaussian-measure}
\begin{split}
d\gamma_{t; q}(y)
& = \frac{\sqrt{1 - q}}{\pi \sqrt{t}} \sin(\theta) (q; q)_\infty \abs{(q e^{2 i \theta}; q)_\infty}^2 dy \\
& = (q; q)_\infty \abs{(q e^{2 i \theta}; q)_\infty}^2 d\gamma_t(\sqrt{1-q} \, y)
\end{split}
\end{equation}
supported on the interval
\[
\left[- \frac{2 \sqrt{t}}{\sqrt{1-q}}, \frac{2 \sqrt{t}}{\sqrt{1-q}}\right].
\]
Here we have used the change of variable \eqref{Eq:q-x-y} and the $q$-Pochhammer symbol
\[
(a_1, \ldots, a_k; q)_\infty = \prod_{j=1}^k \prod_{i=0}^\infty (1 - a_j q^i).
\]
According to Corollary~3.10 of \cite{BKSQGauss}, the $q$-Brownian motion is a Markov process, and moreover the $q$-Hermite polynomials are martingale polynomials with respect to it. Here the (Rogers) continuous $q$-Hermite polynomials
\[
H_n(y, t; q) = t^{n/2} H_n(x/\sqrt{t}; q)
\]
are the monic orthogonal polynomials with respect to the measure \eqref{Eq:q-Gaussian-measure},
\[
\int_{\mf{R}} H_n(y, t; q) H_k(y, t; q) \,d\gamma_{t; q}(y) = \delta_{n=k} [n]_q! t^n.
\]
They also satisfy the three-term recursion relation
\[
y H_n(y, t; q) = H_{n+1} (y, t; q) + [n]_q t H_{n-1}(y, t; q),
\]
where $[n]_q = 1 + q + \ldots + q^{n-1}$.

\begin{Lemma}
The transition operators
\[
\mc{K}_{s, t; q} f(x) = \int_{\mf{R}} f(y) \mc{K}_{s, t; q}(x, dy)
\]
of the $q$-Brownian motion are
\[
\begin{split}
\mc{K}_{s,t; q}(x, dy)
& = \frac{\sqrt{1 - q}}{\pi \sqrt{t}} (q; q)_\infty \sin(\theta) \abs{(q e^{2 i \theta};q)_\infty}^2 \frac{(s/t; q)_\infty}{\abs{(\sqrt{s/t} e^{i(\phi + \theta)}, \sqrt{s/t} e^{i(\phi - \theta)}; q)_\infty}^2} dy, \\
& = (q; q)_\infty \abs{(q e^{2 i \theta};q)_\infty}^2 \frac{(s/t; q)_\infty}{\abs{(\sqrt{s/t} e^{i(\phi + \theta)}, \sqrt{s/t} e^{i(\phi - \theta)}; q)_\infty}^2} d\gamma_t(\sqrt{1 - q} \, y),
\end{split}
\]
where
\begin{equation}
\label{Eq:q-x-y}
x = \frac{2 \sqrt{s}}{\sqrt{1-q}} \cos(\phi), \quad y = \frac{2 \sqrt{t}}{\sqrt{1-q}} \cos(\theta), \quad \phi, \theta \in [0, \pi].
\end{equation}
\end{Lemma}

\begin{proof}
Using the martingale property
\begin{equation}
\label{Eq:Hermite-Martingale}
(\mc{K}_{s,t; q} H_n(y, t; q))(x) = H_n(x, s; q)
\end{equation}
and the orthogonality and density of $q$-Hermite polynomials,
\[
\begin{split}
\mc{K}_{s, t; q}(x, dy)
& = \sum_{n=0}^\infty \frac{1}{[n]_q! t^n} H_n(x, s; q) H_n(y, t; q) \,d\gamma_{t; q}(y) \\
& = \sum_{n=0}^\infty \frac{(s/t)^{n/2}}{[n]_q!} H_n(x/\sqrt{s}; q) H_n(y/\sqrt{t}; q) \,d\gamma_{t; q}(y).
\end{split}
\]
The result now follows from the $q$-Mehler formula
\begin{equation}
\label{Eq:q-Mehler}
\sum_{n=0}^\infty \frac{r^n}{[n]_q!} H_n(x; q) H_n(y; q)
= \frac{(r^2; q)_\infty}{\abs{(r e^{i(\phi + \theta)}, r e^{i(\phi - \theta)}; q)_\infty}^2}.
\end{equation}
See Theorem~4.6 of \cite{BKSQGauss} for more details.
\end{proof}

\begin{Remark}
According to \cite{Donati-Martin}, the It\^{o} product formula for the $q$-Brownian motion has the form
\begin{equation}
\label{Eq:Ito-product}
\begin{split}
\left(\int_0^\infty U(t) \sharp dX(t) \right) \left(\int_0^\infty V(t) \sharp dX(t) \right)
&= \int_0^\infty A(t) d X(t) \left(B(t) \int_0^t V(s) \sharp dX(s) \right) \\
&\quad + \int_0^\infty \left(\int_0^t U(s) \sharp dX(s) \right) C(t) dX(t) D(t) \\
&\quad + \int_0^\infty A(t) \Gamma_q \Bigl[B(t) C(t) \Bigr] D(t) dt,
\end{split}
\end{equation}
where $U = A \otimes B$, $V = C \otimes D$ are adapted biprocesses satisfying a technical condition. Here $\Gamma_q$ is a certain completely positive map on the von Neumann algebra $W^\ast(\set{X(t), t \geq 0})$. In this paper we are only interested in the action of this map on the von Neumann algebra generated by a single operator $X(t)$. This algebra is commutative and isomorphic to
\[
L^\infty\left[- \frac{2 \sqrt{t}}{\sqrt{1-q}}, \frac{2 \sqrt{t}}{\sqrt{1-q}}\right].
\]
On this algebra, the map is determined by the property that
\begin{equation}
\label{Eq:Multiplier}
\Gamma_{t;q}[H_n(x, t; q)] = q^n H_n(x, t; q) = H_n(q x, q^2 t; q),
\end{equation}
that is, it is a multiplier for the $q$-Hermite polynomials. Comparing with equation~\eqref{Eq:Hermite-Martingale}, we see that
 \begin{equation}
\label{Eq:Gamma-K}
\Gamma_{t; q}(x, dy) = \mc{K}_{q^2 t, t; q}(q x, dy).
\end{equation}
In particular, $\Gamma_{t; q}$ is an integral operator
\[
\begin{split}
\Gamma_{t; q}(x, dy)
& = \frac{\sqrt{1 - q}}{\pi \sqrt{t}} \sin(\theta) (q^2; q)_\infty (q; q)_\infty \frac{\abs{(q e^{2 i \theta};q)_\infty}^2}{\abs{(q e^{i(\phi + \theta)}, q e^{i(\phi - \theta)}; q)_\infty}^2} \,dy \\
& = (q^2; q)_\infty (q; q)_\infty \frac{\abs{(q e^{2 i \theta};q)_\infty}^2}{\abs{(q e^{i(\phi + \theta)}, q e^{i(\phi - \theta)}; q)_\infty}^2} d\gamma_t(\sqrt{1 - q} \, y).
\end{split}
\]
\end{Remark}

\begin{Prop}[Functional It\^{o} formula]
\label{Prop:Functional-Ito}
Let $f$ be a polynomial. Then
\begin{equation}
\label{Eq:Functional-Ito}
f(X(t))
= \int_0^t (\partial f)(X(s), X(s)) \,\sharp\, dX(s)
+ \int_0^t (\Delta_{s; q} f)(X(s)) ds.
\end{equation}
Here
\[
\Delta_{s; q} f (x)
= \int_{\mf{R}} \left( \partial_x \frac{f(x) - f(y)}{x - y} \right) \Gamma_{s; q}(x, dy)
= \int_{\mf{R}} (\partial_x \partial f) (x,y) \Gamma_{s; q}(x, dy).
\]
\end{Prop}

\begin{proof}
First we show that all the terms in the functional It\^{o} formula satisfy the technical condition of Theorem~3.2 from \cite{Donati-Martin}. All the integrands are polynomials in $X(s)$. So it suffices to show all the properties for the process $\set{X(t)}$. It is clearly adapted and uniformly bounded on the interval $[0,t]$. Now let
\[
\mc{I} = \set{0 = t_0 < t_1 < \ldots < t_n = t}
\]
be a subdivision of $[0,t]$, and $\delta(\mc{I})$ be the length of the largest interval in this subdivision. Let
\[
X^{\mc{I}}(s) = \sum_{i=0}^{n-1} X(t_i) \chf{[t_i, t_{i+1})}(s).
\]
Then
\[
\int_0^t \norm{X(s) - X^{\mc{I}}(s)}_\infty^2 ds
= \sum_{i=0}^{n-1} \int_{t_i}^{t_{i+1}} \norm{X(s) - X(t_i)}_\infty^2 ds.
\]
But $\norm{X(s) - X(t_i)}_\infty^2 = \norm{X(s - t_i)}_\infty^2 = \frac{4}{1-q} (s - t_i)$. Therefore the preceding sum is
\[
\sum_{i=0}^{n-1} \int_{t_i}^{t_{i+1}} \frac{4}{1-q} (s - t_i) ds
= \frac{2}{1-q} \sum_{i=0}^{n-1} (t_{i+1} - t_i)^2
\leq \frac{2}{1-q} \delta(\mc{I}) \rightarrow 0
\]
as $\delta(\mc{I}) \rightarrow 0$.

The rest of the proof proceeds by induction. Assuming formula \eqref{Eq:Functional-Ito} for $f$, and using the It{\^o} product formula \eqref{Eq:Ito-product}, we get
\[
\begin{split}
f(X(t)) X(t)
&= \int_0^t (\partial f) (X(s), X(s)) \,\sharp (I \otimes X(s)) \,\sharp \,dX(s) + \int_0^t f(X(s)) \,dX(s) \\
&\quad + \int_0^t [(I \otimes \Gamma_q) (\partial f)](X(s)) \,ds + \int_0^t \Delta_{s; q}(f)[X(s)] X(s) \,ds \\
&= \int_0^t (\partial f) (X(s), X(s)) \,\sharp (I \otimes X(s)) \,\sharp \,dX(s) + \int_0^t f(X(s)) \,dX(s) \\
&\quad + \int_0^t [(I \otimes \Gamma_q) (\partial f)](X(s)) \,ds + \int_0^t [(I \otimes \Gamma_q)(\partial_x \partial f)] (X(s)) X(s) \,ds.
\end{split}
\]
The result now follows for $x f(x)$ by observing that
\[
(\partial f)(x,y) y + f(x) = (\partial (x f)) (x,y).
\]
and
\[
(\partial f)(x,y) + \partial_x (\partial f)(x,y) x
= \partial_x(x (\partial f)(x,y))
= \partial_x \Bigl(\partial(x f)(x,y) - f(y) \Bigr)
= \partial_x \partial(x f)(x,y).
\]
\end{proof}

\begin{Cor}
\label{Corollary:q-generator}
On the domain $\mc{P}$ of polynomials, the generators of the $q$-Brownian motion are
\[
A_t f(x) = \Delta_{t; q} f(x) = \int (\partial_x \partial f)(x, y) \Gamma_{t; q}(x, dy),
\]
More explicitly,
\[
\begin{split}
& A_t f(x)
= \int \left(\partial_x \frac{f(x) - f(y)}{x-y} \right) \frac{(q^2; q)_\infty}{\abs{(q e^{i(\phi + \theta)}, q e^{i(\phi - \theta)}; q)_\infty}^2} \;d\gamma_{t; q}(y) \\
&\quad = \int \left(\partial_x \frac{f(x) - f(y)}{x-y} \right)  (q^2; q)_\infty (q; q)_\infty \frac{\abs{(q e^{2 i \theta};q)_\infty}^2}{\abs{(q e^{i(\phi + \theta)}, q e^{i(\phi - \theta)}; q)_\infty}^2} \,d\gamma_t(\sqrt{1-q} y)
\end{split}
\]
with the change of variables \eqref{Eq:q-x-y}.
\end{Cor}

\begin{proof}
It follows from the It\^{o} product formula in Proposition~\ref{Prop:Functional-Ito} that for polynomial $f$,
\[
f(X(t)) - \int_0^t \Delta_{s; q} f(X(s)) \,ds
\]
is a martingale. Therefore by Lemma~\ref{Lemma:Generator-martingale}, $\Delta_{t,q}$ is the generator of the process at time $t$. Note that since the support of $\gamma_{t; q}$ is infinite, polynomials are determined by their values on it. The explicit formula follows.
\end{proof}

\begin{Thm}
\label{Thm:q-BM-generator}
The operator $\Delta_{t; q}$ described in Corollary~\ref{Corollary:q-generator} is the generator of the $q$-Brownian motion at time $t$ on the domain $W^\infty \subset C_0(\mf{R})$.
\end{Thm}

\begin{proof}
First, using the beginning of the proof of Proposition~\ref{Prop:Lp-estimates},
\[
\abs{\Delta_{t; q} f(x)}
= \abs{\int_{\mf{R}} (\partial_x \partial f)(x, y) \Gamma_{t; q}(x, dy)}
\leq \norm{f''}_\infty \abs{\int_{\mf{R}} \Gamma_{t; q}(x, dy)}
\leq \norm{f''},
\]
where in the last step we used equation~\eqref{Eq:Multiplier} for $n=0$. It is also clear that $\mc{K}_{s, t; q}$ is a contraction on $L^\infty(\mf{R}, \,dx)$. By a standard argument, polynomials are dense with respect to the $W^\infty$ norm in $C\left[- \frac{2 \sqrt{t}}{\sqrt{1-q}}, \frac{2 \sqrt{t}}{\sqrt{1-q}}\right]$. Finally, the strong continuity of $\mc{K}_{s, t; q}$ and $\Gamma_{t; q}$ on polynomials follows from the martingale property of the $q$-Hermite polynomials, and formula~\eqref{Eq:Multiplier}. It remains to apply Proposition~\ref{Prop:Basic}.
\end{proof}

\begin{Remark}
Setting $q=0$ in the formula in Corollary~\ref{Corollary:q-generator}, we get
\[
A_t f(x)
= \int \left(\partial_x \frac{f(x) - f(y)}{x-y} \right) \,d\gamma_t(y),
\]
as expected. On the other hand, setting $q=1$ in formula \eqref{Eq:Multiplier} we see that $\Gamma_{t; 1}(x, dy) = \delta(x - y) \,dy$ is the identity operator. So in this case,
\[
\begin{split}
A_t f(x) 
& = \int \left(\partial_x \frac{f(x) - f(y)}{x-y} \right) \delta(x-y) \,dy \\
& = \int \left( \frac{f(y) - f(x) - f'(x) (y-x)}{(y-x)^2} \right) \delta(x-y) \,dy 
= \frac{1}{2} f''(x),
\end{split}
\]
again as expected.
\end{Remark}

\section{Two-state free Brownian motions}
\label{Section:Two}

In \cite{Ans-Two-Brownian}, we considered Brownian motions in the context of two-state free probability theory $(\mc{A}, \mf{E}, E)$. A priori, any process with two-state freely independent increments whose $E$-distributions are a free convolution semigroup $\set{\nu_t}$ and whose $\mf{E}$-distributions satisfy
\begin{equation*}
\label{Eq:Two-state-Gaussian}
\mc{J}[\mu_t] = \nu_t, \quad \mu_t[x] = 0, \mu_t[x^2] = t
\end{equation*}
can be considered a two-state free Brownian motion. We proved, however, that if we require $\mf{E}$ to be a faithful normal state and $E$ be normal, then $\nu_t$ has to be the semicircular distribution with possibly non-zero mean $\alpha t$ and variance $t$. Such a process is not Markov (in fact $\mf{E}$ is not tracial, and $\mf{E}$-preserving conditional expectations do not exist), however its classical version is. We now construct generators of these processes.

\begin{Prop}
The generator of the two-state free Brownian motion $\set{X(t)}$ with parameter $\alpha$ at time $t$ is
\[
\alpha (\partial_x - L_{\mu_t}) + \partial_x L_{\nu_t},
\]
on the domain $W^{\infty}$.
\end{Prop}

\begin{proof}
This result was proved in \cite{Ans-Two-Brownian} for polynomial $f$. Also,
\[
\norm{(\alpha (\partial_x - L_{\mu_t}) + \partial_x L_{\nu_t}) f} \leq 2 \abs{\alpha} \norm{f'}_\infty + \norm{f''}_\infty \leq 2 \abs{\alpha} \norm{f}_\infty + (2 \abs{\alpha} + 1) \norm{f''}_\infty.
\]
Since the measures $\mu_t, \nu_t$ are all uniformly compactly supported, the full result follows as in Theorem~\ref{Thm:q-BM-generator}.
\end{proof}

\begin{Remark}[It{\^o} formula]
By the same methods as in \cite{BiaSpeBrownian} and \cite{AnsIto}, for sufficiently nice $f$,
\begin{equation}
\label{Eq:Ito}
f(X(t)) = f(0) + \int_0^t \partial f(X(s)) \sharp \,dX(s) + \int_0^t (\partial_x \otimes E) \partial f(X(s)) \,ds.
\end{equation}
Using Lemma~2.1 of \cite{BLS96} and the observation that the process $\set{X(t)}$ is $\mf{E}$-centered (see Remark~6 of \cite{Ans-Two-Brownian} for more details), we see that
\[
\Exp{f(X(t))}
= f(0) + \int_0^t \Exp{ \Bigl( \alpha \partial_x - \alpha (1 \otimes \mf{E}) \partial + (\partial_x \otimes E) \partial \Bigr) f(X(s))} \,ds.
\]
This result is consistent with the generator formula in the preceding proposition.
\end{Remark}


\def\cprime{$'$}
\providecommand{\bysame}{\leavevmode\hbox to3em{\hrulefill}\thinspace}
\providecommand{\MR}{\relax\ifhmode\unskip\space\fi MR }
\providecommand{\MRhref}[2]{%
  \href{http://www.ams.org/mathscinet-getitem?mr=#1}{#2}
}
\providecommand{\href}[2]{#2}

\end{document}